\documentclass[11pt,a4paper]{amsart}
\usepackage[utf8]{inputenc}
\usepackage[T1]{fontenc}
\usepackage{amsmath,amssymb,amsthm,mathtools}
\usepackage[margin=2.6cm]{geometry}
\usepackage{xcolor}
\usepackage[authoryear,round]{natbib}
\usepackage{hyperref}
\usepackage{enumitem}
\usepackage[english]{babel}
\usepackage{pgfplots}
\usepgfplotslibrary{groupplots}
\pgfplotsset{compat=1.17}
\hypersetup{colorlinks=true,linkcolor=blue!50!black,citecolor=blue!50!black,urlcolor=blue!50!black}

\numberwithin{equation}{section}
\allowdisplaybreaks

\theoremstyle{plain}
\newtheorem{thm}{Theorem}[section]
\newtheorem{prop}[thm]{Proposition}
\newtheorem{lem}[thm]{Lemma}
\newtheorem{cor}[thm]{Corollary}
\theoremstyle{definition}
\newtheorem{defn}[thm]{Definition}
\newtheorem{assum}[thm]{Assumption}
\newtheorem{exmp}[thm]{Example}
\theoremstyle{remark}

\newenvironment{pf}{\begin{proof}}{\end{proof}}

\newcommand{\figwidth}{9.5cm}

\newcommand{\X}{X}
\newcommand{\U}{U}
\newcommand{\UU}{\mathbb{U}}
\newcommand{\Useq}{\mathcal{U}}
\newcommand{\Uad}{\mathcal{U}_F}
\newcommand{\Mcl}{M_{\mathrm{cl}}}
\newcommand{\Scl}{S_{\mathrm{cl}}}
\newcommand{\PF}{P_F}

\newcommand{\eps}{\varepsilon}
\DeclarePairedDelimiter{\norm}{\lVert}{\rVert}
\newcommand{\ip}[2]{\langle #1,#2\rangle}
\DeclareMathOperator{\ess}{ess}
\newcommand{\zet}{\zeta}

\title[Cost controllability and predictive control of monotone systems]{Cost controllability and predictive control of input-constrained monotone systems on Hilbert spaces}
\author[Till Preuster]{Till Preuster\textsuperscript{1}}
\author[Karl Worthmann]{Karl Worthmann\textsuperscript{2}}
\author[Lars Gr\"une]{Lars Gr\"une\textsuperscript{3}}
\author[Manuel Schaller]{Manuel Schaller\textsuperscript{1}}

\thanks{\textsuperscript{1}Faculty of Mathematics, Chemnitz University of Technology,
	Germany (\texttt{\{till.preuster,manuel.schaller\}@math.tu-chemnitz.de}).}
\thanks{\textsuperscript{2}Institute of Mathematics, Technische Universit\"at Ilmenau,
	Germany (\texttt{karl.worthmann@tu-ilmenau.de}).}
\thanks{\textsuperscript{3}Mathematical Institute, University of Bayreuth, Germany
	(\texttt{lars.gruene@uni-bayreuth.de}).}
\thanks{Till Preuster and Manuel Schaller acknowledge funding by the Deutsche Forschungsgemeinschaft (DFG, German Research Foundation) – Project-ID 531152215 – CRC 1701.}
\date{\today}

\begin{document}

	\begin{abstract}
		Cost controllability, relating controllability and a performance measure like the stage cost in model predictive control (MPC),
		is a key tool in the stability and suboptimality analysis of closed-loop systems.
		We establish
		cost controllability
		for maximal monotone
		control systems subject to input constraints, a class encompassing finite- and
		infinite-dimensional port-Hamiltonian systems with colocated input and output.
		To this end, we employ a recently proposed saturated feedback
		to derive state-dependent bounds on the closed-loop cost.
		Our main result is that
		cost controllability holds globally with semiglobal exponential rate.
		We first establish this result in continuous time and then extend it to sampled-data systems.
		Whereas a qualitative detectability-like
		property is sufficient for asymptotic convergence of the controller, cost controllability requires quantitative rates. We provide two verifiable conditions ensuring cost
		controllability with explicit exponential rates: coercivity of the sum of
		system operator and output injection, and exponential detectability of the linear part.
		Moreover, under the coercivity condition, the resulting MPC scheme without
		terminal conditions renders the equilibrium globally asymptotically stable with semiglobal exponential rate already for
		prediction horizon two.
		We illustrate our
		results by a nonlinear finite-dimensional system, a semilinear heat equation and a
		nonlinearly damped wave equation.
	\end{abstract}

	\maketitle


	\section{Introduction}\label{sec:intro}

	Model predictive control (MPC) computes a feedback law by iteratively solving finite-horizon
	optimal control problems and applying the first portion of the optimal control function, see,
	e.g., \citep{GruenePannek,RawlingsMayneDiehl}. Asymptotic stability of the MPC closed loop may be
	enforced by terminal costs and constraints, which, however, are often hard to construct, in
	particular for nonlinear partial differential equations. Without (stabilizing) terminal
	conditions, the stability and performance analysis of MPC \citep{grimm2005model,Gruene09} usually
	rests on a controllability assumption in terms of the stage cost, called \emph{cost controllability}.
	Loosely speaking, it requires, for every initial state, an admissible control such that the stage
	cost along the resulting trajectory decays according to a fixed $\mathcal{KL}_0$-bound depending
	only on the initial cost. For bounds linear in the initial cost, the minimal stabilizing horizon
	can be computed by a linear program \citep{GruenePannekSeehafer10,Worthmann11}, and this analysis was extended
	to continuous time in \citet{RebleAllgoewer12}.

	Verifying cost controllability in applications can be challenging, in particular under input
	constraints, and considerable work has been devoted to particular system classes.
	\citet{CoronGrueneWorthmann20} derived sufficient conditions in terms of homogeneous
	approximations, \citet{WorthmannMehrez16} established cost controllability for nonholonomic
	mobile robots by explicit open-loop maneuvers, and \citet{MuellerWorthmann17} showed that the
	design of the stage cost is decisive, as quadratic costs may fail to asymptotically stabilize the
	MPC closed loop, whereas suitably designed costs succeed \citep{CoronGrueneWorthmann20,rosenfelder2023model}. In infinite
	dimensions, cost controllability and MPC without terminal conditions were studied for heat and wave
	equations in \citet{AltmuellerGruene12,AltmuellerGrueneWorthmann10}, with emphasis on stage costs
	allowing for small stabilizing horizons. Receding horizon stabilization of infinite-dimensional
	systems was considered in \citet{ItoKunisch02} and established for Burgers and wave equations in
	\citet{AzmiKunisch16,AzmiKunisch18}. A second line of research connects MPC with
	passivity and dissipativity. Closed-loop stability may be enforced by a passivity-based constraint
	\citep{RaffEbenbauerAllgoewer07}, dissipativity is central to economic MPC, see, e.g.,
	\citep{AngeliAmritRawlings12}, and for port-Hamiltonian, hence dissipative, systems, optimal
	control problems with the supplied energy as cost exhibit turnpike phenomena
	\citep{SchallerPhilipp21}.

	Maximal monotone operators provide a unifying framework for nonlinear and
	infinite-dimen\-sional dissipative dynamics, e.g., in heat conduction, wave propagation or
	port-Hamilto\-nian networks, cf.\ \citep{Barbu,CamlibelSchaft}. The stabilization of such systems
	under a priori bounds on the control goes back to \citet{Slemrod89}. Cone-bounded and saturated
	damping injection was analyzed for the wave equation in \citet{PrieurTarbouriech16}, for the
	Korteweg-de Vries equation in \citet{MarxCerpa17} and for general m-dissipative nonlinear systems
	in \citet{MarxAndrieuPrieur17}. \citet{Vanspranghe25} consider projected dynamic feedback of
	impedance passive systems. For maximal monotone systems with colocated input and output,
	\citet{PGS26} recently proposed a saturated feedback, that is, the output damping injection
	projected onto the admissible control values. It respects arbitrary closed convex constraint sets
	and, under a detectability-like condition, renders the controlled equilibrium asymptotically
	stable.

	\noindent\textbf{Contribution and outline.} We show that maximal monotone control systems under input
	constraints are globally cost controllable with semiglobal exponential rate, such that MPC stability theory directly becomes
	applicable. Our contributions are threefold. (i) Using the saturated feedback as admissible
	control, we prove cost controllability in continuous time (Theorems \ref{thm:coercive} and~\ref{thm:detect}) and for the sampled-data system
	(Corollary~\ref{cor:discrete}), with explicit semiglobal exponential rates under coercivity or
	exponential detectability. (ii) Under coercivity of $M+BB^*$, MPC without terminal conditions renders the equilibrium
	globally asymptotically stable with semiglobal exponential rate
	already for the shortest horizon $N=2$ and every sampling period (Theorem~\ref{thm:N2}). (iii) The results apply to nonlinear partial differential
	equations, as illustrated by a semilinear heat equation and a nonlinearly damped wave equation
	(Examples~\ref{ex:heat} and~\ref{ex:wave}). After recalling the system class and the saturated
	feedback in Section~\ref{sec:prelim}, we reduce cost controllability to a state decay estimate in
	Section~\ref{sec:cc}, verify this estimate in Section~\ref{sec:decay}, derive the MPC
	guarantees in Section~\ref{sec:mpc} and conclude in Section~\ref{sec:conclusion}.

	\noindent\textbf{Notation.} Throughout, $\X$ and $\U$ are real Hilbert spaces with inner
	products $\ip{\cdot}{\cdot}$ and $\ip{\cdot}{\cdot}_\U$ and norms $\norm{\cdot}$ and $\norm{\cdot}_\U$, and
	$\mathcal{L}(\U,\X)$ denotes the space of bounded linear operators.
	The open ball of radius $r$ around $u$ is $B_r(u)$. A
	continuous function $\rho:\mathbb{R}_{\ge0}\to\mathbb{R}_{\ge0}$ is of class
	$\mathcal{K}_\infty$ if it is strictly increasing, unbounded and satisfies $\rho(0)=0$. A continuous function
	$\beta:\mathbb{R}_{\ge0}\times\mathbb{R}_{\ge0}\to\mathbb{R}_{\ge0}$ is of class
	$\mathcal{KL}_0$ if $\beta(r,t)\to0$ for $t\to\infty$ for each $r>0$ and if, for each
	$t\ge0$, either $\beta(\cdot,t)\in\mathcal{K}_\infty$ or $\beta(\cdot,t)\equiv0$. For $\bar\rho>0$, we write
	$\overline{D(M)}_{\bar\rho}\coloneqq\{x\in\overline{D(M)}:\norm{x-x^\star}\le\bar\rho\}$.

	\section{Maximal monotone control systems and saturated feedback}\label{sec:prelim}

	We briefly recall maximal monotone control systems and the saturated feedback, which serves as an auxiliary tool for our main results.

	\begin{defn}
		A possibly unbounded and nonlinear operator $M:\X\supset D(M)\to\X$ is called \emph{monotone} if
		\begin{equation*}
			\langle M(x)-M(z),x-z\rangle_\X \geq 0, \quad \text{for all } x,z\in D(M).
		\end{equation*}
		It is called \emph{maximal monotone} if, in addition, it admits no proper monotone extension or, equivalently, if $\lambda I+M$ is onto for some (and hence all) $\lambda>0$; cf.\ \citep{Barbu} for more background on maximal monotone operators and their associated nonlinear semigroups.
	\end{defn}

	For $M:\X\supset D(M)\to\X$ maximal monotone and $B\in\mathcal{L}(\U,\X)$,
	we consider the \emph{maximal monotone control system} $(M,B)$ given by
	\begin{align}
		\label{eq:sys}
		\begin{split}
			\tfrac{\mathrm{d}}{\mathrm{d}t}x(t)=-M(x(t))+Bu(t), \quad y(t)=B^*x(t),
		\end{split}
	\end{align}
	where, for initial conditions $x(0)=x_0\in\overline{D(M)}$ and integrable inhomogeneities,
	solutions are understood in the mild sense, cf.\ \citep[Thm.~4.5]{Barbu}. The class~\eqref{eq:sys}
	encompasses linear dissipative systems, nonlinear gradient-type flows and finite- as well as
	infinite-dimensional port-Hamiltonian models.
	We refer the reader to
	\citet{CamlibelSchaft} for the close relation between port-Hamiltonian systems and monotonicity. Consider a closed convex
	set of control constraints~$F\subset\U$, and a \textit{controlled equilibrium} $(x^\star,u^\star)\in D(M)\times F$,
	i.e. $-M(x^\star)+Bu^\star=0$.
	The saturated feedback proposed in \citet{PGS26} reads
	\begin{equation}\label{eq:fb}
		\kappa(x)\coloneqq \PF\big(u^\star-B^*(x-x^\star)\big)\in F ,
	\end{equation}
	where $\PF$ denotes the $\U$-orthogonal projection onto $F$. For unconstrained problems,
	that is, $F=\U$, the projection is the identity and \eqref{eq:fb} reduces to the
	colocated damping injection $u=u^\star-B^*(x-x^\star)$ which is the standard stabilizing feedback for
	port-Hamiltonian systems, cf.\ \citep{ortega2002stabilization}.
	By \citet[Lems.~3.2, 3.3]{PGS26}, the corresponding closed-loop
	operator $\Mcl(x)\coloneqq M(x)-B\kappa(x)$ is maximal monotone with
	$\Mcl(x^\star)=0$ such that $-\Mcl$ generates a semigroup
	$(\Scl(t))_{t\ge0}$ of nonlinear nonexpansive mappings on $\overline{D(M)}$. Since $\Scl(t)x^\star=x^\star$ for all $t\ge0$, the nonexpansiveness and the semigroup property yield, for all $x_0\in\overline{D(M)}$ and $0\le s\le t$,
	\begin{equation}\label{eq:mono}
		\norm{\Scl(t)x_0-x^\star}\le\norm{\Scl(s)x_0-x^\star}\le\norm{x_0-x^\star}.
	\end{equation}
	The stability analysis of the closed-loop rests on the following detectability-like property.

	\begin{defn}[Vanishing output vanishing state property; VOVS]\label{def:vovs}
		A maximal monotone system
		$ \dot x=-M(x),\ y=B^*x $,
		with semigroup $(S(t))_{t\ge0}$ generated by
		$-M$ has the \emph{vanishing output vanishing state} property at $x^\star$ if,
		for all $x_0\in\overline{D(M)}$, $\lim_{t\to\infty}B^*(S(t)x_0-x^\star)=0$ implies
		$\lim_{t\to\infty}S(t)x_0=x^\star$.
	\end{defn}

	By \citep[Thm.~3.7]{PGS26}, the VOVS property yields global asymptotic stability of the closed loop.

	\begin{thm}[Asymptotic stability]\label{thm:PGS}
		Let $u^\star\in\operatorname{int}F$. If the closed-loop system $\dot x=-M_\mathrm{cl}(x), y=B^*x$, fulfills
		the VOVS property at $x^\star$, then $\Scl(t)x_0\to x^\star$ as $t\to\infty$ for all $x_0\in\overline{D(M)}$. Precisely, $x^\star$ is globally asymptotically stable, for the closed loop resulting
		from
		the feedback laws~\eqref{eq:fb}, i.e., plugging in $u(t)=\kappa(x(t))$ into \eqref{eq:sys}.
	\end{thm}

	The VOVS property is a purely qualitative detectability notion, imposing neither a rate nor
	any uniformity in the initial value. Quantitative counterparts are provided by output-to-state
	stability \citep{SontagWang97,KrichmanSontagWang01,ChenMironchenkoWirth25}, and for linear
	systems, exponential detectability implies the VOVS property of both the open loop and the closed loop subject to \eqref{eq:fb}, cf.\ \citet[Props.~3.8, 3.9]{PGS26}.

	\section{Cost controllability}\label{sec:cc}
	In this part, we formulate cost controllability for the system \eqref{eq:sys} in
	continuous time and prove the first main result of this work: along the saturated feedback
	\eqref{eq:fb}, the control part of the stage cost is dominated by the state, such that
	cost controllability reduces to a uniform decay estimate for the closed-loop state.

	As all mild solutions of \eqref{eq:sys} evolve in $\overline{D(M)}\subset\X$, the set $\overline{D(M)}$ will serve
	as the state space throughout. We define the set of admissible controls
	$\Uad\coloneqq \{u\in L^2_{\mathrm{loc}}([0,\infty);\U):u(t)\in F \text{ for a.e.\ } t\ge0\}$.
	For $x_0\in\overline{D(M)}$ and $u\in\Uad$ we denote by $x_u(\cdot\,;x_0)$ the unique mild solution of
	\eqref{eq:sys}, which is well defined due to $Bu\in L^1_{\mathrm{loc}}([0,\infty);\X)$, cf.\
	\citep[Thm.~4.5]{Barbu}, and takes values in $\overline{D(M)}$. As stage cost we choose the quadratic function $\ell: X \times U \to \mathbb{R}$ with
	\begin{subequations}\label{eq:ell}
		\begin{align}
			\ell(x,u) &\coloneqq \norm{x-x^\star}^2 + \lambda\norm{u-u^\star}_\U^2, \label{eq:ell-a}\\
			\ell^*(x) &\coloneqq \inf_{u\in F}\ell(x,u) = \norm{x-x^\star}^2, \label{eq:ell-b}
		\end{align}
	\end{subequations}
	with $\lambda > 0$ where the infimum is attained at $u=u^\star\in F$. The following  is
	a central assumption underlying the continuous-time MPC
	analysis, see \citep{RebleAllgoewer12,CoronGrueneWorthmann20}.

	\begin{defn}[Cost controllability]\label{def:cc}
		Let $\beta\in\mathcal{KL}_0$ be such that $\beta(r,\cdot) \in L^1(0, \infty)$ for every $r\ge0$. The system \eqref{eq:sys} is
		\emph{cost controllable} at $(x^\star,u^\star)\in D(M)\times F$ with
		rate $\beta$ if for each $x_0\in\overline{D(M)}$ there is
		$u_{x_0}\in\Uad$ such that
		\begin{equation}\label{eq:cc}
			\ell\big(x_{u_{x_0}}(t;x_0),u_{x_0}(t)\big)\le\beta\big(\ell^*(x_0),t\big),
			\quad t\ge0 .
		\end{equation}
	\end{defn}
	Under cost controllability, integrating \eqref{eq:cc} yields
	a bound on the finite-horizon value function
	$V_t(x_0) \coloneqq \inf_{u\in\Uad}\int_0^t\ell(x_u(s;x_0),u(s))\,\mathrm ds$,
	that is,
	\begin{equation}\label{eq:Vbound}
		V_t(x_0)
		\le\int_0^\infty\beta(\ell^*(x_0),s)\,\mathrm ds<\infty
		\quad\text{for all } t\ge0 .
	\end{equation}
	In particular, if $\beta(r,t)=r\,b(t)$, $V_t(x_0)\le \int_0^t b(s)\,\mathrm ds\,\ell^*(x_0)$.

	Cost controllability is the key hypothesis in the stability analysis of MPC without terminal conditions.
	The main idea of this work is to verify it by employing the feedback \eqref{eq:fb} as the admissible control in
	Definition~\ref{def:cc}.
	A central tool is \textit{firm nonexpansivity} of the projection $\PF$, see \citep{BausComb2011}, i.e., for all $v,w\in\U$,
	\begin{equation}\label{eq:firm}
		\ip{\PF(v)-\PF(w)}{v-w}_\U\ge\norm{\PF(v)-\PF(w)}_\U^2.
	\end{equation}

	\begin{lem}\label{lem:adm}
		Let $x_0\in\overline{D(M)}$, let $x(t)\coloneqq \Scl(t)x_0$ be the closed-loop trajectory and for $t\geq 0$,
		\begin{equation}\label{eq:ucl}
			u_{x_0}(t)\coloneqq \kappa(x(t))=\PF\big(u^\star-B^*(x(t)-x^\star)\big).
		\end{equation}
		Then $u_{x_0}$ is continuous and bounded with $u_{x_0}(t)\in F$ for all $t\ge0$; in
		particular $u_{x_0}\in\Uad$. Moreover, $u_{x_0}$ reproduces the closed loop, that is,
		$x_{u_{x_0}}(\cdot\,;x_0)=x$. If further
		$x_0\in D(M)$, then $u_{x_0}$ is Lipschitz continuous and $x$
		is a strong solution of \eqref{eq:sys} with $u=u_{x_0}$.
	\end{lem}

	\begin{pf}
		The projection $\PF$ is firmly nonexpansive by \eqref{eq:firm} and maps into $F$, and
		$t\mapsto x(t)$ is continuous, hence $u_{x_0}$ is continuous with values in $F$. The
		stability estimate \eqref{eq:mono} and $\PF(u^\star)=u^\star$ yield
		$\norm{u_{x_0}(t)-u^\star}\le\norm{B^*}\norm{x_0-x^\star}$, hence $u_{x_0}$ is bounded.
		If $x_0\in D(M)$, then $t\mapsto x(t)$ is Lipschitz continuous with
		$\norm{\dot x(t)}\le\norm{\Mcl(x_0)}$ for a.e.\ $t\geq 0$, \citep[Thm.~4.4]{Barbu}, hence so is $u_{x_0}$, and
		$Bu_{x_0}\in W^{1,\infty}_{\mathrm{loc}}([0,\infty);\X)$, such that $x_{u_{x_0}}(\cdot\,;x_0)$ is a
		strong solution by \citet[Thm.~4.5]{Barbu}. Moreover, since $-\Mcl(x)=-M(x)+B\kappa(x)$ by
		definition of $\Mcl$, the closed-loop trajectory $x$ is a strong solution of
		$\dot x=-M(x)+Bu_{x_0}(t)$ with $x(0)=x_0$, and uniqueness of mild solutions yields
		$x_{u_{x_0}}(\cdot\,;x_0)=x$. For general $x_0\in\overline{D(M)}$, this identity follows by
		density, as $\Scl(t)$ and $\PF$ are nonexpansive and mild solutions depend continuously on
		initial values and inhomogeneities, cf.\ \citep{Barbu}.
	\end{pf}

	The key estimate for the proposed control is the following dissipation inequality.
	\begin{lem}[Dissipation inequality]\label{lem:diss}
		Let $x_0\in\overline{D(M)}$, $x(t)=\Scl(t)x_0$ and let $u_{x_0}$ be defined by \eqref{eq:ucl}. Then, for
		all $0\le s\le t<\infty$,
		\begin{equation*}
			\tfrac12\norm{x(t)-x^\star}^2 + \!\int_s^t\norm{u_{x_0}(\tau)-u^\star}_\U^2\,\mathrm{d}\tau
			\le\tfrac12\norm{x(s)-x^\star}^2 .
		\end{equation*}
		In particular,
		$\int_s^\infty\norm{u_{x_0}(\tau)-u^\star}^2\,\mathrm{d}\tau\le\tfrac12\norm{x(s)-x^\star}^2$, and
		$t\mapsto\norm{x(t)-x^\star}$ is non-increasing.
	\end{lem}

	\begin{pf}
		We first prove the claim for $x_0\in D(M)$. In this case $x$ is the strong solution, which
		is Lipschitz continuous on $[0,\infty)$ and hence differentiable almost everywhere.
		Abbreviating $\delta x(t) \coloneqq x(t)-x^\star$, $\delta u(t) \coloneqq u_{x_0}(t)-u^\star$ and
		$v(t)\coloneqq u^\star-B^*\delta x(t)$, we have $\delta u(t)=\PF(v(t))-\PF(u^\star)$ by
		\eqref{eq:ucl} and $u^\star\in F$. Using $\Mcl(x)=M(x)-B\kappa(x)$ and $\Mcl(x^\star)=0$, we
		obtain for a.e.\ $t$:
		\begin{equation}\label{eq:dissipation-core}
			\begin{aligned}
				\tfrac{\mathrm d}{\mathrm{d}t}\tfrac12\norm{\delta x(t)}^2
				&=-\ip{M(x(t))-M(x^\star)}{\delta x(t)} +\ip{\delta u(t)}{B^*\delta x(t)}_\U ,
			\end{aligned}
		\end{equation}
		where the first term is nonpositive by monotonicity. Further, the firm nonexpansiveness \eqref{eq:firm} yields
		\begin{align}\label{eq:key-sector}
			\begin{split}
				&\ip{\delta u(t)}{B^*\delta x(t)}_\U=-\ip{\delta u(t)}{v(t)-u^\star}_\U \\
				&\hspace{3cm} = \!-\langle P_F(v(t)) \!-\! P_F(u^\star), v(t) \!-\!u^\star\rangle_U\! \le \!-\norm{\delta u(t)}_\U^2 ,
			\end{split}
		\end{align}
		hence $\frac{\mathrm d}{\mathrm{d}t}\tfrac12\norm{\delta x(t)}^2\le-\norm{\delta u(t)}^2$ for
		a.e.\ $t\ge0$, and integration over $[s,t]$ gives the claimed estimate. For $x_0\in\overline{D(M)}$, choose $x_0^k\in D(M)$ with $x_0^k\to x_0$. By nonexpansivity of $\Scl(t)$ and $\PF$, the states and control deviations converge pointwise, $x^k(t)\to x(t)$ and $\delta u^k(t)\to\delta u(t)$, with $\norm{\delta u^k(t)}\le\norm{B^*}\sup_k\norm{x_0^k-x^\star}$. Dominated convergence on $[s,t]$ yields the estimate for $x_0$, and letting $t\to\infty$ gives the $L^2(s,\infty)$-bound.
	\end{pf}

	We now show that cost controllability follows from a decay estimate on the state.
	\begin{prop}\label{prop:reduction}
		Suppose that there is $\tilde\beta\in\mathcal{KL}_0$ with
		$\tilde\beta(r,\cdot)\in L^1(0,\infty)$ for every $r\ge0$ such that, for all
		$x_0\in\overline{D(M)}$ and $t\ge0$,
		\begin{equation}\label{eq:KLstate}
			\norm{\Scl(t)x_0-x^\star}^2\le\tilde\beta\big(\norm{x_0-x^\star}^2,\,t\big) .
		\end{equation}
		Then the system \eqref{eq:sys} is cost controllable at $(x^\star,u^\star)$ with
		the admissible control \eqref{eq:ucl} and the rate
		\begin{equation}\label{eq:beta-pw}
			\beta(r,t)\coloneqq \big(1+\lambda\norm{B^*}^2\big)\,\tilde\beta(r,t) ,
		\end{equation}
		which belongs to $\mathcal{KL}_0$ and satisfies $\beta(r,\cdot)\in L^1(0,\infty)$ for every $r\ge0$. In particular,
		\eqref{eq:Vbound} holds. Conversely, if~\eqref{eq:sys} is cost controllable at
		$(x^\star,u^\star)$ with rate $\beta$, then for every $x_0\in\overline{D(M)}$
		there is $u\in\Uad$ such that
		$\norm{x_u(t;x_0)-x^\star}^2\le\beta(\norm{x_0-x^\star}^2,t)$ for all $t\ge0$.
		In other words, cost controllability implies uniform asymptotic controllability to $x^\star$ by $F$-valued controls with integrable rate, and, conversely, a uniform decay estimate with integrable rate along the feedback \eqref{eq:fb} implies cost controllability.
	\end{prop}

	\begin{pf}
		Let $x_0\in\overline{D(M)}$ and let $u_{x_0}$ be given by \eqref{eq:ucl}, which
		is admissible by Lemma~\ref{lem:adm} and generates the trajectory
		$x(t)=x_{u_{x_0}}(t;x_0)=\Scl(t)x_0$. Since $\PF$ is nonexpansive by
		\eqref{eq:firm} and $\PF(u^\star)=u^\star$, we have
		$\norm{u_{x_0}(t)-u^\star}_\U\le\norm{B^*\delta x(t)}_\U\le\norm{B^*}\norm{\delta x(t)}$
		with $\delta x(t)=x(t)-x^\star$, hence
		$\ell(x(t),u_{x_0}(t))\le(1+\lambda\norm{B^*}^2)\norm{\delta x(t)}^2
		\le(1+\lambda\norm{B^*}^2)\,\tilde\beta(\ell^*(x_0),t)$
		by \eqref{eq:KLstate} and $\ell^*(x_0)=\norm{x_0-x^\star}^2$. Since
		multiplication by a positive constant preserves $\mathcal{KL}_0$ and
		integrability, this proves cost controllability with the rate \eqref{eq:beta-pw}.
		The converse follows from $\ell(x,u)\ge\norm{x-x^\star}^2$ and the definition
		of $\ell^*$.
	\end{pf}

	\section{Verification of the decay estimate~\eqref{eq:KLstate}}\label{sec:decay}
	By Proposition~\ref{prop:reduction}, the uniform decay estimate \eqref{eq:KLstate} along the saturated feedback \eqref{eq:fb} implies cost controllability. Theorem~\ref{thm:PGS} only gives convergence for each initial value separately, without a rate. We thus provide two quantitative conditions yielding an integrable $\tilde\beta\in\mathcal{KL}_0$. Both use a local sector condition, which controls the effect of the input constraints. Together with this sector condition, semiglobal coercivity (Subsection~\ref{ssec:coercive}) or exponential detectability of the linear part (Subsection~\ref{ssec:detect}) yield global decay with semiglobal exponential rate.

	By \eqref{eq:firm}, the control deviation $\kappa(x)-u^\star$ is never larger than the output deviation $B^*(x-x^\star)$. For the decay estimates we additionally need a lower bound of the same order, which depends on the geometry of $F$ and holds if, e.g., $u^\star\in\operatorname{int}F$. It ensures that a vanishing control
	deviation forces a vanishing output and hence, under additional detectability-like assumptions, a vanishing state.
	\begin{assum}[Local sector condition]\label{ass:S}
		There is a continuous non-increasing function $\theta:[0,\infty)\to(0,1]$ such that for all $\rho \geq 0$ and $x\in\overline{D(M)}$ with $\norm{x-x^\star}\le\rho$,
		\[
		\norm[\big]{\PF\big(u^\star-B^*(x-x^\star)\big)-u^\star}_\U\ge\theta(\rho)\,\norm{B^*(x-x^\star)}_\U .
		\]
	\end{assum}
	The following result shows that Assumption~\ref{ass:S} quantifies the interiority condition
	$u^\star\in\operatorname{int}F$ that is common in control-constrained stabilization.

	\begin{lem}\label{lem:sector}
		Let $F\subset\U$ be closed and convex, $u^\star\in F$ and $r>0$ with $B_r(u^\star)\subset F$,
		that is, $u^\star\in\operatorname{int}F$. Then the following holds.
		\begin{enumerate}[label=(\roman*)]
			\item For every $v\in\U$, either $\PF(v)=v$ or $\norm{\PF(v)-u^\star}_\U\ge r$; consequently, for all $v\in\U$
			\begin{equation}\label{eq:sector}
				\norm{\PF(v)-u^\star}_\U\ \ge\ \min\{r,\norm{v-u^\star}_\U\}.
			\end{equation}
			\item Assumption~\ref{ass:S} holds with $\theta(\rho)=\min\{1,\,r/(\norm{B^*}\rho)\}$ for
			$\rho>0$ and $\theta(0)\coloneqq 1$.
		\end{enumerate}
	\end{lem}

	\begin{pf}
		We first show (i). Let $q\coloneqq \PF(v)$ and assume $q\ne v$ and $\norm{q-u^\star}<r$. Then
		$q\in B_r(u^\star)$, which is open and contained in $F$, hence
		$z\coloneqq q+\eps\frac{v-q}{\norm{v-q}}\in F$ for sufficiently small $\eps>0$. The variational
		characterization of the projection, $\ip{v-q}{z-q}_\U\le0$ for all $z\in F$, then yields
		$0\ge\ip{v-q}{\eps\frac{v-q}{\norm{v-q}}}=\eps\norm{v-q}>0$, a contradiction. Inequality
		\eqref{eq:sector} follows since $\PF(v)=v$ gives $\norm{\PF(v)-u^\star}=\norm{v-u^\star}$.
		Concerning (ii), let $\norm{x-x^\star}\le\rho$ and put $v\coloneqq u^\star-B^*(x-x^\star)$, such that
		$\eta\coloneqq \norm{v-u^\star}=\norm{B^*(x-x^\star)}\le\norm{B^*}\rho$. If $\eta\le r$, then
		\eqref{eq:sector} gives $\norm{\PF(v)-u^\star}\ge\eta$; if $\eta>r$, then \eqref{eq:sector}
		gives
		$\norm{\PF(v)-u^\star}\ge r=\frac{r}{\norm{B^*}\rho}\norm{B^*}\rho\ge\frac{r}{\norm{B^*}\rho}\,\eta$.
	\end{pf}
	For constraints defined by pointwise bounds, which are ubiquitous in problems governed by PDEs, the interiority condition of Lemma~\ref{lem:sector} is void.
	The next result provides the
	appropriate substitute.

	\begin{lem}[Pointwise box constraints]\label{lem:box}
		Let $\Omega_c\subset\mathbb{R}^d$ be measurable with positive measure, $\U=L^2(\Omega_c)$,
		$F=\{u\in\U: u_{\min}\le u\le u_{\max}\text{ a.e.}\}$, and suppose that
		$$\rho_0\coloneqq \ess\inf_{\xi\in\Omega_c}\min\{u^\star(\xi)-u_{\min},\,u_{\max}-u^\star(\xi)\}>0.$$
		Then, for every $R>0$, all $x\in\overline{D(M)}$ with $\norm{B^*(x-x^\star)}_{L^\infty(\Omega_c)}\le R$ satisfy the estimate of Assumption~\ref{ass:S} with $\theta(\rho)$ replaced by $\min\{1,\rho_0/R\}$.
	\end{lem}

	\begin{pf}
		The projection acts pointwise, that is,
		$\PF(v)(\xi)=\min\{u_{\max},\max\{u_{\min},v(\xi)\}\}$. Applying
		Lemma~\ref{lem:sector}(i) in $\U=\mathbb{R}$ with the interval
		$F_\xi=[u_{\min},u_{\max}]$, which contains
		$[u^\star(\xi)-\rho_0,u^\star(\xi)+\rho_0]$, gives,
		with $\eta\coloneqq B^*(x-x^\star)$,
		$\big|\PF(u^\star-\eta)(\xi)-u^\star(\xi)\big|\ge\min\{\rho_0,|\eta(\xi)|\}\ge\min\Big\{1,\tfrac{\rho_0}{R}\Big\}\,|\eta(\xi)|$,
		where the second inequality distinguishes the cases $|\eta(\xi)|\le\rho_0$ and
		$|\eta(\xi)|>\rho_0$ and uses $|\eta(\xi)|\le R$ a.e. Squaring and integrating over $\Omega_c$ proves the claim.
	\end{pf}

	\subsection{Exponential decay under quantitative VOVS conditions}\label{sec:sharp}

	As the VOVS property imposes neither a rate nor uniformity in the initial value, cf.\
	Section~\ref{sec:prelim}, we strengthen it in two ways, each giving \eqref{eq:KLstate} with explicit exponential rates, namely by coercivity of $M+BB^*$
	(Subsection~\ref{ssec:coercive}) and by exponential detectability
	of the linear part (Subsection~\ref{ssec:detect}).

	\subsubsection{Coercive systems}\label{ssec:coercive}
	We first assume that $M+BB^*$ is coercive at $x^\star$ on bounded sets, in which case a
	Gr\"onwall argument yields global decay with semiglobal exponential rate.
	\begin{assum}[Semiglobal coercivity]\label{ass:D}
		There is a continuous non-increasing function $\zet:[0,\infty)\to(0,\infty)$ such that for
		all $\rho\ge0$ and $x\in D(M)$ with $\norm{x-x^\star}\le\rho$,
		\begin{equation}\label{eq:D}
			\begin{aligned}
				&\ip{M(x)-M(x^\star)}{x-x^\star}+\norm{B^*(x-x^\star)}_\U^2\\
				&\qquad\ge\zet(\rho)\,\norm{x-x^\star}^2 .
			\end{aligned}
		\end{equation}
		If $\zet$ can be chosen constant, we call the system globally coercive. Further, let
		$\zet_M\ge0$ be a constant such that
		\begin{equation}\label{eq:DM}
			\ip{M(x)-M(x^\star)}{x-x^\star}\ \ge\ \zet_M\,\norm{x-x^\star}^2
		\end{equation}
		for all $x\in D(M)$.
		Due to monotonicity of $M$, $\zet_M=0$ is always admissible, and since
		\eqref{eq:DM} implies \eqref{eq:D} with $\zet_M$ in place of $\zet(\rho)$, we will assume w.l.o.g.\ $\zet\ge\zet_M$.
	\end{assum}

	We note that \eqref{eq:D} implies the VOVS property as in
	Definition~\ref{def:vovs} for the closed-loop system $(\Mcl,B)$, such that Assumption~\ref{ass:D} may be viewed as a
	quantitative strengthening of the VOVS property. The constant $\zet_M$ isolates the portion of the decay that is provided by the plant itself. In the rate below, only the remaining, output-injection part of \eqref{eq:D} is affected by the saturation of the input. The following result shows that Assumption~\ref{ass:D} yields exponential decay.

	\begin{thm}[Exponential decay under coercivity]\label{thm:coercive}
		Let Assumptions~\ref{ass:S} and \ref{ass:D} hold and set
		\begin{equation}\label{eq:mudef}
			\mu(\rho)\coloneqq \big(1-\theta(\rho)^2\big)\,\zet_M+\theta(\rho)^2\zet(\rho)\ \ge\ \zet_M .
		\end{equation}
		Then, for every $x_0\in\overline{D(M)}$ and with
		$\rho\coloneqq \norm{x_0-x^\star}$,
		\begin{equation}\label{eq:expdecay}
			\norm{\Scl(t)x_0-x^\star}^2\ \le\ e^{-2\mu(\rho)\,t}\,\rho^2\qquad\text{for all } t\ge0 .
		\end{equation}
		Further,
		\eqref{eq:sys} is cost controllable at $(x^\star,u^\star)$
		with rate
		\begin{equation}\label{eq:beta-coercive}
			\beta(r,t)=\big(1+\lambda\norm{B^*}^2\big)\,e^{-2\mu(\sqrt r)\,t}\,r .
		\end{equation}
		In particular, on every ball $\overline{D(M)}_{\bar\rho}$, the decay is exponential with the uniform rate
		$\mu(\bar\rho)$.
	\end{thm}

	\begin{pf}
		We first prove \eqref{eq:expdecay} for $x_0\in D(M)$ and abbreviate $x(t)=\Scl(t)x_0$,
		$\delta x(t)=x(t)-x^\star$, $\delta y(t)=B^*\delta x(t)$ and
		$\delta u(t)=\PF(u^\star-\delta y(t))-u^\star$. Due to \eqref{eq:mono} we have
		$\norm{\delta x(t)}\le\rho$ for all $t\ge0$, such that Assumptions~\ref{ass:S} and \ref{ass:D}
		apply with this fixed $\rho$ along the whole trajectory. Using \eqref{eq:dissipation-core},
		\eqref{eq:key-sector}, the bound $\norm{\delta u}^2\ge\theta(\rho)^2\norm{\delta y}^2$ from
		Assumption~\ref{ass:S} and splitting the monotonicity term according to
		$1=(1-\theta(\rho)^2)+\theta(\rho)^2$, we obtain for a.e.\ $t$, omitting the argument $t$,
		\begin{align*}
			\tfrac{\mathrm d}{\mathrm{d}t}\tfrac12\norm{\delta x}^2&\le-\ip{M(x)-M(x^\star)}{\delta x}-\norm{\delta u}^2\\
			&\le-\big(1-\theta(\rho)^2\big)\ip{M(x)-M(x^\star)}{\delta x} - \theta(\rho)^2\big(\ip{M(x)-M(x^\star)}{\delta x}+\norm{\delta y}^2\big)\\
			&\le-\mu(\rho)\norm{\delta x}^2 ,
		\end{align*}
		where the last inequality follows from \eqref{eq:DM} and \eqref{eq:D}. Gr\"onwall's inequality gives
		\eqref{eq:expdecay}; for $x_0\in\overline{D(M)}$, choose $x_0^k\in D(M)$ with $x_0^k\to x_0$. Since $\Scl(t)$ is nonexpansive and $\mu$ is continuous, both sides of \eqref{eq:expdecay} converge, which yields the estimate for $x_0$. Hence, \eqref{eq:KLstate} holds with $\tilde\beta(r,t)\coloneqq e^{-2\mu(\sqrt r)t}r$, and Proposition~\ref{prop:reduction} yields cost controllability at $(x^\star,u^\star)$ with the rate \eqref{eq:beta-coercive}, provided that $\tilde\beta\in\mathcal{KL}_0$. To show this,
		we first note that $\mu(\rho)=\zet_M+\theta(\rho)^2(\zet(\rho)-\zet_M)>0$ is
		continuous and non-increasing, being, due to $\zet\ge\zet_M$, the sum of a constant and
		a product of two nonnegative, continuous, non-increasing functions. Hence
		$(r,t)\mapsto e^{-2\mu(\sqrt r)t}$ is continuous and, for fixed $t\ge0$, non-decreasing in
		$r$ with values in $(0,1]$, such that $r\mapsto e^{-2\mu(\sqrt r)t}r$
		belongs to $\mathcal{K}_\infty$. For
		fixed $r>0$, $e^{-2\mu(\sqrt r)t}r\to0$ as $t\to\infty$, hence $\tilde\beta\in\mathcal{KL}_0$.
		Finally, $\mu(\rho)\ge\mu(\bar\rho)$ for $x_0\in\overline{D(M)}_{\bar\rho}$, as
		$\mu$ is non-increasing, such that \eqref{eq:expdecay} gives the uniform rate on balls.
	\end{pf}
	The rate \eqref{eq:beta-coercive} is valid for every initial value without smallness
	condition, although the controls are constrained to $F$, as monotonicity of $M$ ensures
	that the bounded control never has to dominate unstable growth. For $\zet_M=0$, however, the
	exponent $\mu(\sqrt r)$ degrades as the initial deviation grows, since far from the equilibrium
	the feedback saturates, such that the decay is global with semiglobal exponential rate, i.e.,
	exponential bounds with uniform constants hold on every ball.

	\begin{exmp}[Nonlinear finite-dimensional system]\label{ex:rates}
		Let $\X=\U=\mathbb{R}^2$ and $M(x)=c\norm{x-x^\star}^2(x-x^\star)+J(x-x^\star)$ with
		$c>0$ and $J=\left(\begin{smallmatrix}0&1\\-1&0\end{smallmatrix}\right)$. Then $M$ is the
		sum of the gradient of the convex function $\frac{c}{4}\norm{x-x^\star}^4$ and a
		skew-symmetric linear map, hence continuous and monotone, and thus maximal monotone.
		Further, let $B=I$, $F=[a,b]^2$ with $a<0<b$ and $u^\star=0$, so that
		$(x^\star,u^\star)$ is a controlled equilibrium. Without control, the state decays only
		polynomially, since $\frac{\mathrm d}{\mathrm dt}\norm{x(t)-x^\star}^2=-2c\norm{x(t)-x^\star}^4$
		yields $\norm{x(t)-x^\star}^2=\rho_0^2/(1+2c\rho_0^2t)$ with $\rho_0=\norm{x_0-x^\star}$.
		With $z=x-x^\star$ and $\ip{Jz}{z}=0$, we have
		\[
		\ip{M(x)-M(x^\star)}{z}+\norm{B^*z}^2=c\norm{z}^4+\norm{z}^2\ge\norm{z}^2,
		\]
		so Assumption~\ref{ass:D} holds with $\zet\equiv1$ and $\zet_M=0$. Since
		$B_{r_F}(u^\star)\subset F$ for $r_F\coloneqq\min\{-a,b\}$ and $\norm{B^*}=1$,
		Lemma~\ref{lem:sector} yields Assumption~\ref{ass:S} with
		$\theta(\rho)=\min\{1,r_F/\rho\}$. Hence, Theorem~\ref{thm:coercive} gives
		\[
		\norm{x(t)-x^\star}^2\le\exp\big(-2\min\{1,r_F^2/\rho_0^2\}\,t\big)\,\rho_0^2,
		\]
		and cost controllability with the integrable rate
		$\beta(r,t)=(1+\lambda)\,e^{-2\min\{1,r_F^2/r\}t}\,r$. The decay is thus exponential
		on every ball, with exponent $2$ for $\rho_0\le r_F$, which deteriorates due to the
		saturation for $\rho_0>r_F$.
	\end{exmp}

	\begin{exmp}[Semilinear heat equation with Neumann boundary conditions]\label{ex:heat}
		We consider $M=A_N+\mathcal N$ on $\X=L^2(\Omega)$ with the Neumann Laplacian $A_Nx=-\Delta x$, $D(A_N)=\{x\in H^2(\Omega):\partial_\nu x=0\}$, on the unit square $\Omega$ and the Nemytskii operator $\mathcal Ny\coloneqq f\circ y$,
		$D(\mathcal N)\coloneqq\{y\in L^2(\Omega): f\circ y\in L^2(\Omega)\}$, of a continuous and nondecreasing $f\colon\mathbb{R}\to\mathbb{R}$ with $f(0)=0$, such that $M$ is maximal monotone with $D(M)=D(A_N)$ by Proposition~\ref{prop:AN+F}. Further, let $\U=L^2(\Omega_c)$, $B=\mathbf 1_{\Omega_c}$ the extension by zero from the control region $\Omega_c$, such that $B^*x=x|_{\Omega_c}$ and $\norm{B}=\norm{B^*}=1$, $F\coloneqq\{u\in\U: T_{\min}\le u\le T_{\max}\text{ a.e.}\}$, and let $(x^\star,u^\star)\in D(M)\times F$ be a controlled equilibrium. By the coercivity of
		$A_N+BB^*$ shown in the proof of \citet[Prop.~4.3]{PGS26} and the monotonicity of $f$, Assumption~\ref{ass:D} holds with $\zet=\min\{\frac{|\Omega_c|}{4(1+2/C_P)},\frac{C_P}{2}\}$, where $C_P$ is the Poincar\'e constant of $\Omega$. If $f$ is strongly monotone, that is, $(f(a)-f(b))(a-b)\ge m(a-b)^2$ for some $m>0$, then \eqref{eq:D} and \eqref{eq:DM} hold with
		$\zet+m$ and $\zet_M=m$. As $F$ has empty $L^2$-interior, Assumption~\ref{ass:S} fails on every ball of $\X$, and we use Lemma~\ref{lem:box} with $\rho_0=\ess\inf\min\{u^\star-T_{\min},T_{\max}-u^\star\}$ instead, e.g., $\rho_0\ge1$ for the data $T_{\min}=-5$, $T_{\max}=7$ and $\norm{u^\star}_{L^\infty}\le4$. Let $R_\infty>0$. For every $x_0\in\overline{D(M)}$ with $\norm{x_0-x^\star}_{L^\infty(\Omega)}\le R_\infty$, Lemma~\ref{lem:heatLinf} yields
		\[
		\norm{B^*(\Scl(t)x_0-x^\star)}_{L^\infty(\Omega_c)}\le\norm{\Scl(t)x_0-x^\star}_{L^\infty(\Omega)}\le R_\infty
		\]
		for all $t\ge0$. Moreover, $x_0^k\coloneqq x^\star+(I+\tfrac1kW)^{-1}(x_0-x^\star)\in D(M)$ with $W$ from the proof of Lemma~\ref{lem:heatLinf} satisfies $\norm{x_0^k-x^\star}_{L^\infty(\Omega)}\le R_\infty$ by that proof and $x_0^k\to x_0$ as $k\to\infty$, cf.\ \citep{Barbu}, such that the same bound holds along $\Scl(t)x_0^k$. As the proof of Theorem~\ref{thm:coercive} uses Assumption~\ref{ass:S} only along the closed-loop trajectories from $x_0^k$, Lemma~\ref{lem:box} and Theorem~\ref{thm:coercive} with $\theta\equiv\min\{1,\rho_0/R_\infty\}$ yield \eqref{eq:expdecay} with $\mu=m+\zet\min\{1,\rho_0^2/R_\infty^2\}$, where $m=0$ for merely monotone $f$, and, due to
		$\norm{B^*}=1$, cost controllability with the rate $\beta(r,t)=(1+\lambda)e^{-2\mu t}r$ for these initial values. Note that this is weaker than cost controllability in the sense of Definition~\ref{def:cc}, since it only covers initial values in the $L^\infty(\Omega)$-ball of radius $R_\infty$ and the rate depends, via $\mu$, on $R_\infty$ and not only on $\ell^*(x_0)$.
	\end{exmp}

	\subsubsection{Exponentially detectable linear part}\label{ssec:detect}

	Coercivity of $M+BB^*$ may fail if $B^*$ observes only part of the state, e.g., for the wave equation in Example~\ref{ex:wave} below, where damping and output act on the velocity only. In this part, we therefore split $M$ into a linear part, which is only required to be exponentially detectable, and a co-coercive nonlinearity, and derive the decay estimate from the variation-of-parameters formula combined with the dissipation inequality of Lemma~\ref{lem:diss}. Precisely, we impose the following.
	\begin{assum}[Exponentially detectable linear part]\label{ass:det}
		Let $M=A+N$, where $A$ is a closed, densely defined, linear maximal
		monotone operator (so that $-A$ generates a linear contraction semigroup),
		and $N:\X\to\X$ is co-coercive, that is, there is $c_N>0$ such that
		\begin{equation}\label{eq:cocoercive}
			\ip{N(x)-N(z)}{x-z}\ \ge\ c_N\norm{N(x)-N(z)}^2
		\end{equation}
		for all $x,z\in\X$. Further, $(A,B)$ is \textit{exponentially detectable}, that is, there exist
		$L\in\mathcal{L}(\U,\X)$, $c\ge1$ and $\omega>0$ such that
		$\norm{T_{-A+LB^*}(t)}\le ce^{-\omega t}$ for all $t\ge0$.
	\end{assum}
	Co-coercive operators are monotone and Lipschitz continuous with constant
	$1/c_N$. The case $N=0$ is admissible with the convention $c_N^{-1/2}\coloneqq0$
	and covers linear systems. Moreover, $M$ is maximal monotone with $D(M)=D(A)$,
	\citep[Cor.~25.5]{BausComb2011}.

	We will now show that this assumption of exponential detectability of the linear part yields global decay with a quantitative, semiglobal exponential rate.

	\begin{thm}[Exponential decay under detectability]\label{thm:detect}
		Let Assumptions~\ref{ass:S} and~\ref{ass:det} be fulfilled. Define, for $\rho\ge0$,
		$g(\rho)\coloneqq \norm{B}+\norm{L}/\theta(\rho)$, $\nu(\rho)\coloneqq \ln 2/(4K(\rho)^2)$ with
		\begin{equation}\label{eq:Krho}
			K(\rho)\coloneqq \tfrac{c}{\sqrt2}\Big(\tfrac{1}{\sqrt\omega}+\tfrac{g(\rho)+c_N^{-1/2}}{\omega}\Big) .
		\end{equation}
		Then, for every $x_0\in\overline{D(M)}$ and for all $t\geq 0$,
		\begin{equation}\label{eq:detect-decay}
			\norm{\Scl(t)x_0-x^\star}\le 2e^{-\nu(\norm{x_0-x^\star})t}\norm{x_0-x^\star}.
		\end{equation}
		Consequently, \eqref{eq:KLstate} holds with
		$\tilde\beta(r,t)=4e^{-2\nu(\sqrt r)t}r$, which belongs to $\mathcal{KL}_0$
		and satisfies $\tilde\beta(r,\cdot)\in L^1(0,\infty)$ for every $r\ge0$, such that \eqref{eq:sys} is cost
		controllable at $(x^\star,u^\star)$ with rate
		$\beta(r,t)=4\big(1+\lambda\norm{B^*}^2\big)e^{-2\nu(\sqrt r)t}r$.
	\end{thm}

	The proof of Theorem~\ref{thm:detect} is given in Appendix~\ref{app:detect}.
	In the linear and unconstrained case, that is, $N=0$ and
	$F=\U$ with $\theta\equiv1$, Theorem~\ref{thm:detect} recovers the well-known exponential stability of $-(A+BB^*)$ for
	exponentially detectable $(A,B)$, here with an explicit rate. For general constraint sets $F$, the input constraint enters the rate only via the replacement of $\norm{L}$ by $\norm{L}/\theta(\rho)$ in $g(\rho)$.

	\begin{exmp}[Nonlinearly damped wave equation]\label{ex:wave}
		Let $\Omega\subset\mathbb{R}^d$ be a bounded domain with smooth boundary and let the open control region $\Omega_c\subset\Omega$ satisfy the geometric control condition (GCC) of \citet{BardosLebeauRauch92}, that is, there is $T_0>0$ such that every ray of geometric optics in $\Omega$, travelling at unit speed and reflected at $\partial\Omega$, enters $\Omega_c$ before time $T_0$. We consider the wave equation with nonlinear damping
		\[
		\partial_{tt}w-\Delta w+\psi(\partial_tw)=\mathbf 1_{\Omega_c}u ,
		\]
		with homogeneous Dirichlet boundary conditions, where $\psi:\mathbb{R}\to\mathbb{R}$ is
		nondecreasing and Lipschitz continuous with constant $L_\psi$ and $\psi(0)=0$, e.g.,
		$\psi=\arctan$. With $z=(w,\partial_tw)\in\X=H_0^1(\Omega)\times L^2(\Omega)$ and
		$\U=L^2(\Omega_c)$, this is \eqref{eq:sys} with $M=A+N$, where $A$ is the skew-adjoint wave
		operator, $N(z)=(0,\psi\circ z_2)$ and $B^*z=z_2|_{\Omega_c}$, such
		that $\norm{B}=\norm{B^*}=1$ and $(x^\star,u^\star)=(0,0)$ is a controlled equilibrium. Since
		$\psi$ is nondecreasing and $L_\psi$-Lipschitz, we have
		$(\psi(a)-\psi(b))(a-b)\ge|\psi(a)-\psi(b)|^2/L_\psi$ for all $a,b\in\mathbb{R}$, such that $N$
		satisfies \eqref{eq:cocoercive} with $c_N=1/L_\psi$. As both the damping and the output act on
		the velocity only, the coercivity condition of Assumption~\ref{ass:D} fails. However, by the GCC, $(A,B)$ is exactly controllable \citep{BardosLebeauRauch92} and hence exponentially
		detectable. For the energy constraint
		$F=\{u\in\U:\norm{u}_\U\le u_{\max}\}$ with $u_{\max}>0$, we have $B_{u_{\max}}(u^\star)\subset F$,
		such that Lemma~\ref{lem:sector} yields Assumption~\ref{ass:S} with
		$\theta(\rho)=\min\{1,u_{\max}/\rho\}$. Hence, Theorem~\ref{thm:detect} yields global
		cost controllability with semiglobal exponential rate determined by \eqref{eq:Krho}, where
		$c_N^{-1/2}=L_\psi^{1/2}$.
	\end{exmp}

	\section{Application to model predictive control}\label{sec:mpc}

	In this part, we derive closed-loop guarantees for MPC without terminal conditions. The works analyzing
	continuous-time MPC that are recalled briefly before Definition~\ref{def:cc} rely on cost
	controllability but treat unconstrained
	finite-dimensional dynamics. For the constrained Hilbert space setting of this work, we instead
	employ the discrete-time framework of \citet{Gruene09}, which is formulated on arbitrary metric
	spaces and accommodates constraints. To this end, consider the discrete-time system
	\begin{equation}\label{eq:dt_system}
		x(n+1)=f(x(n),u(n))
	\end{equation}
	on the state space $\overline{D(M)}$, equipped with the metric of $\X$, with control values in
	a set $\UU$, control sequences $\Useq\coloneqq \{u:\mathbb{N}_0\to\UU\}$, solutions $x_u(n)$ with
	$x_u(0)=x_0$ and a running cost $l$ with $l^*(x)\coloneqq\min_{u\in\UU}l(x,u)$. Given a horizon $N\in\mathbb{N}$ and $x_0\in\overline{D(M)}$, we consider
	\begin{equation}\label{eq:ocp}
		V_N(x_0)\coloneqq\inf_{u\in\UU^N}\sum_{n=0}^{N-1}l\big(x_u(n),u(n)\big),
	\end{equation}
	where $u=(u(0),\dots,u(N-1))$. At each time $n$, MPC solves \eqref{eq:ocp} with $x_0=x(n)$ and applies the first element $\mu_N(x(n))\coloneqq\hat u(0)$ of an optimal control sequence $(\hat u(0),\dots,\hat u(N-1))$, that is, $x(n+1)=f(x(n),\mu_N(x(n)))$.
	For simplicity, we assume that a minimizer exists, which has to be verified for specific nonlinear partial differential equations, see \citet[Ch.~5]{Troeltzsch10}.

	The discrete-time counterpart of Definition~\ref{def:cc} plays a central role in the stability analysis of this MPC scheme.

	\begin{defn}[Cost controllability in discrete time, {\citealp[cf.][Ass.~3.1]{Gruene09}}]\label{def:ccd}
		Let $\beta\in\mathcal{KL}_0$ be such that $\sum_{n=0}^\infty\beta(r,n)<\infty$ for every $r\ge0$. The discrete-time system \eqref{eq:dt_system} is \emph{cost controllable} with
		rate $\beta$ if for each $x_0\in\overline{D(M)}$ there exists a control
		$u_{x_0}\in\Useq$ such that
		\[
		l\big(x_{u_{x_0}}(n),u_{x_0}(n)\big)\le\beta\big(l^*(x_0),n\big)\qquad\text{for all }n\in\mathbb{N}_0 .
		\]
	\end{defn}

	Under cost controllability, summing the estimate in Definition~\ref{def:ccd} yields, analogously to \eqref{eq:Vbound},
	\begin{equation}\label{eq:VNbound}
		V_N(x_0)\le\sum_{n=0}^{\infty}\beta\big(l^*(x_0),n\big)<\infty\quad\text{for all }N\in\mathbb{N}.
	\end{equation}
	For exponential cost controllability, that is, $\beta(r,n)=C\sigma^nr$ with $C\ge1$ and
	$\sigma\in(0,1)$,
	\citet[Thm.~5.2]{Gruene09} yields asymptotic stability of the MPC closed loop
	for all sufficiently large horizons $N$.

	\subsection{Sampled-data formulation}\label{ssec:sampled}

	To cast \eqref{eq:sys} in this framework, we fix a sampling period $T>0$ and choose the set of
	control values
	\begin{equation}\label{eq:XU}
		\UU\coloneqq \big\{u\in L^2((0,T);\U)\,:\,u(t)\in F \text{ for a.e.\ } t\big\} ,
	\end{equation}
	the transition map $f(x,u)\coloneqq \varphi(T;x,u)$, where $\varphi(\cdot\,;x,u)$ denotes the
	mild solution of \eqref{eq:sys} on $[0,T]$ with initial value $x$ and control $u\in\UU$, and
	the running cost
	\begin{equation}\label{eq:l}
		l(x,u)\coloneqq \norm{x-x^\star}_X^2+\lambda\int_0^T\norm{u(t)-u^\star}_\U^2\,\mathrm{d}t .
	\end{equation}
	Unlike the continuous-time cost $\ell$ from \eqref{eq:ell-a}, which is evaluated pointwise in time, the sampled-data cost $l$ takes a control function on $(0,T)$ as argument, evaluates the state only at the beginning of the sampling interval and integrates the control cost over it. Since $u^\star\in F$, the constant control $u\equiv u^\star$ belongs to $\UU$ and minimizes $l(x,\cdot)$, such that $l$ and $\ell$ have the same minimal value, that is, for all $x\in\overline{D(M)}$,
	\begin{align}\label{eq:lstar}
		\norm{x-x^\star}^2 = l^*(x)=\min_{u\in\UU}l(x,u) =\ell^*(x).
	\end{align}
	Here, a control value $u\in\UU$ is a control function on one sampling interval, which is covered by the framework of \citet{Gruene09}, as it allows for general sets of control values, cf.\ also the receding horizon frameworks for PDEs of \citet{ItoKunisch02,AzmiKunisch16,AzmiKunisch18}. In this way, every admissible control, in
	particular the MPC feedback, satisfies the input constraint, and the restrictions
	$u_{x_0}(nT+\cdot)|_{[0,T)}$ of the feedback \eqref{eq:ucl} are admissible control values by
	Lemma~\ref{lem:adm}.

	\subsection{Sampled-data cost controllability with long horizons}\label{ssec:long}

	The stability analysis of MPC without terminal conditions rests on Definition~\ref{def:ccd}, which follows
	from \eqref{eq:KLstate}.

	\begin{cor}[Sampled-data cost controllability]\label{cor:discrete}
		Suppose that \eqref{eq:KLstate} holds with $\tilde\beta\in\mathcal{KL}_0$ such that $\sum_{n=0}^\infty\tilde\beta(r,nT)<\infty$ for every $r\ge0$. Then the sampled-data problem with \eqref{eq:XU} and \eqref{eq:l} is cost controllable in the sense of Definition~\ref{def:ccd} with the control
		\eqref{eq:ucl} and the rate $\beta(r,n)=(1+\tfrac\lambda2)\,\tilde\beta(r,nT)$.	For $\bar\rho>0$ and $r\le\bar\rho^2$, we have $\beta(r,n)\le C\sigma^nr$ for all $n\in\mathbb{N}_0$, i.e., exponential cost controllability on $\overline{D(M)}_{\bar\rho}$, where, under Assumptions~\ref{ass:S} and~\ref{ass:D}, $C=1+\tfrac\lambda2,\sigma=e^{-2\mu(\bar\rho)T}$ and, under Assumptions~\ref{ass:S} and~\ref{ass:det}, $C=4\big(1+\tfrac\lambda2\big), \sigma=e^{-2\nu(\bar\rho)T}$.
	\end{cor}

	\begin{pf}
		By Lemma~\ref{lem:adm}, $u_{x_0}(nT+\cdot)|_{[0,T)}\in\UU$ for all
		$n\in\mathbb{N}_0$. Thus, $u_{x_0}$ defines an admissible control
		sequence with $x_{u_{x_0}}(n)=\Scl(nT)x_0$. Hence by Lemma~\ref{lem:diss} with
		$s=nT$, $t=(n+1)T$,
		\[
		\int_{nT}^{(n+1)T}\norm{u_{x_0}(\tau)-u^\star}_\U^2\,\mathrm{d}\tau
		\le\tfrac12\norm{\Scl(nT)x_0-x^\star}^2 .
		\]
		By \eqref{eq:l}, \eqref{eq:KLstate} and \eqref{eq:lstar},
		\begin{align*}
			l\big(x_{u_{x_0}}(n),u_{x_0}(n)\big) & \le(1+\tfrac\lambda2)\norm{\Scl(nT)x_0-x^\star}^2 \le(1+\tfrac\lambda2)\,\tilde\beta\big(l^*(x_0),nT\big),
		\end{align*}
		and $\beta$ inherits the $\mathcal{KL}_0$ properties and the summability from $\tilde\beta$.
		For the exponential rates, Theorems~\ref{thm:coercive} and~\ref{thm:detect}
		give $\tilde\beta(r,t)=e^{-2\mu(\sqrt r)t}r$ and
		$\tilde\beta(r,t)=4e^{-2\nu(\sqrt r)t}r$, respectively. As $\mu$ and $\nu$ are
		non-increasing, $\mu(\sqrt r)\ge\mu(\bar\rho)$ and
		$\nu(\sqrt r)\ge\nu(\bar\rho)$ for $r\le\bar\rho^2$, which yields the claim.
	\end{pf}

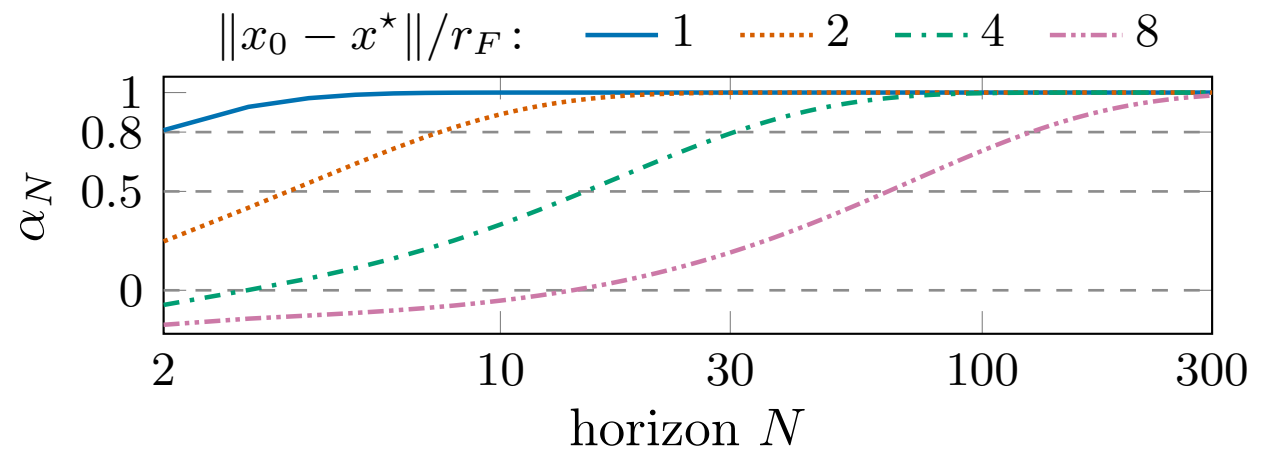
\begin{figure}[t]
		\centering
		\begin{minipage}{\figwidth}\centering
\definecolor{c0}{RGB}{0,114,178}\definecolor{c1}{RGB}{213,94,0}\definecolor{c2}{RGB}{0,158,115}\definecolor{c3}{RGB}{204,121,167}\definecolor{c4}{RGB}{60,60,60}
\begin{tikzpicture}
\begin{axis}[
	scale only axis, width=0.73\linewidth, height=1.7cm,
	xmode=log, xmin=2, xmax=300, ymin=-0.22, ymax=1.08, clip=false,
	xtick={2,10,30,100,300}, xticklabels={$2$,$10$,$30$,$100$,$300$},
	ytick={0,0.5,0.8,1},
	xlabel={horizon $N$}, ylabel={$\alpha_N$}, xlabel style={yshift=3pt}, ylabel style={yshift=-2pt},
	tick label style={font=\footnotesize}, label style={font=\small},
	legend style={at={(0.5,1.03)}, anchor=south, legend columns=-1, draw=none,
	              font=\footnotesize, inner sep=0pt,
	              /tikz/every even column/.append style={column sep=6pt}},
	legend image post style={xscale=0.8},
	every axis plot/.append style={line width=0.85pt},
]
\foreach \a in {0,0.5,0.8}{\addplot[black!45, dashed, line width=0.5pt, forget plot, domain=2:300, samples=2] {\a};}
\addlegendimage{empty legend}\addlegendentry{$\norm{x_0-x^\star}/r_F\colon$}
\addplot[c0, solid] coordinates {(2,0.8097) (3,0.9276) (4,0.9716) (5,0.9888) (6,0.9956) (7,0.9982) (8,0.9993) (9,0.9997) (10,0.9999) (11,1.0000) (12,1.0000) (13,1.0000) (14,1.0000) (15,1.0000) (16,1.0000) (17,1.0000) (18,1.0000) (19,1.0000) (20,1.0000) (21,1.0000) (22,1.0000) (23,1.0000) (24,1.0000) (25,1.0000) (26,1.0000) (27,1.0000) (28,1.0000) (29,1.0000) (30,1.0000) (31,1.0000) (32,1.0000) (33,1.0000) (34,1.0000) (35,1.0000) (36,1.0000) (37,1.0000) (38,1.0000) (39,1.0000) (40,1.0000) (41,1.0000) (42,1.0000) (43,1.0000) (44,1.0000) (45,1.0000) (46,1.0000) (47,1.0000) (48,1.0000) (49,1.0000) (50,1.0000) (51,1.0000) (52,1.0000) (53,1.0000) (54,1.0000) (55,1.0000) (56,1.0000) (57,1.0000) (58,1.0000) (59,1.0000) (60,1.0000) (61,1.0000) (62,1.0000) (64,1.0000) (65,1.0000) (66,1.0000) (68,1.0000) (69,1.0000) (71,1.0000) (72,1.0000) (74,1.0000) (75,1.0000) (77,1.0000) (78,1.0000) (80,1.0000) (81,1.0000) (83,1.0000) (85,1.0000) (87,1.0000) (88,1.0000) (90,1.0000) (92,1.0000) (94,1.0000) (96,1.0000) (98,1.0000) (100,1.0000) (102,1.0000) (104,1.0000) (106,1.0000) (108,1.0000) (111,1.0000) (113,1.0000) (115,1.0000) (118,1.0000) (120,1.0000) (122,1.0000) (125,1.0000) (128,1.0000) (130,1.0000) (133,1.0000) (136,1.0000) (138,1.0000) (141,1.0000) (144,1.0000) (147,1.0000) (150,1.0000) (153,1.0000) (156,1.0000) (160,1.0000) (163,1.0000) (166,1.0000) (170,1.0000) (173,1.0000) (177,1.0000) (180,1.0000) (184,1.0000) (188,1.0000) (192,1.0000) (196,1.0000) (200,1.0000) (204,1.0000) (208,1.0000) (212,1.0000) (217,1.0000) (221,1.0000) (226,1.0000) (230,1.0000) (235,1.0000) (240,1.0000) (245,1.0000) (250,1.0000) (255,1.0000) (260,1.0000) (265,1.0000) (271,1.0000) (277,1.0000) (282,1.0000) (288,1.0000) (294,1.0000) (300,1.0000)};
\addlegendentry{$1$}
\addplot[c1, densely dotted] coordinates {(2,0.2470) (3,0.4175) (4,0.5426) (5,0.6397) (6,0.7160) (7,0.7760) (8,0.8233) (9,0.8606) (10,0.8900) (11,0.9132) (12,0.9315) (13,0.9460) (14,0.9574) (15,0.9663) (16,0.9734) (17,0.9790) (18,0.9835) (19,0.9869) (20,0.9897) (21,0.9919) (22,0.9936) (23,0.9949) (24,0.9960) (25,0.9968) (26,0.9975) (27,0.9980) (28,0.9984) (29,0.9988) (30,0.9990) (31,0.9992) (32,0.9994) (33,0.9995) (34,0.9996) (35,0.9997) (36,0.9998) (37,0.9998) (38,0.9999) (39,0.9999) (40,0.9999) (41,0.9999) (42,0.9999) (43,1.0000) (44,1.0000) (45,1.0000) (46,1.0000) (47,1.0000) (48,1.0000) (49,1.0000) (50,1.0000) (51,1.0000) (52,1.0000) (53,1.0000) (54,1.0000) (55,1.0000) (56,1.0000) (57,1.0000) (58,1.0000) (59,1.0000) (60,1.0000) (61,1.0000) (62,1.0000) (64,1.0000) (65,1.0000) (66,1.0000) (68,1.0000) (69,1.0000) (71,1.0000) (72,1.0000) (74,1.0000) (75,1.0000) (77,1.0000) (78,1.0000) (80,1.0000) (81,1.0000) (83,1.0000) (85,1.0000) (87,1.0000) (88,1.0000) (90,1.0000) (92,1.0000) (94,1.0000) (96,1.0000) (98,1.0000) (100,1.0000) (102,1.0000) (104,1.0000) (106,1.0000) (108,1.0000) (111,1.0000) (113,1.0000) (115,1.0000) (118,1.0000) (120,1.0000) (122,1.0000) (125,1.0000) (128,1.0000) (130,1.0000) (133,1.0000) (136,1.0000) (138,1.0000) (141,1.0000) (144,1.0000) (147,1.0000) (150,1.0000) (153,1.0000) (156,1.0000) (160,1.0000) (163,1.0000) (166,1.0000) (170,1.0000) (173,1.0000) (177,1.0000) (180,1.0000) (184,1.0000) (188,1.0000) (192,1.0000) (196,1.0000) (200,1.0000) (204,1.0000) (208,1.0000) (212,1.0000) (217,1.0000) (221,1.0000) (226,1.0000) (230,1.0000) (235,1.0000) (240,1.0000) (245,1.0000) (250,1.0000) (255,1.0000) (260,1.0000) (265,1.0000) (271,1.0000) (277,1.0000) (282,1.0000) (288,1.0000) (294,1.0000) (300,1.0000)};
\addlegendentry{$2$}
\addplot[c2, dashdotted] coordinates {(2,-0.0741) (3,0.0014) (4,0.0594) (5,0.1120) (6,0.1611) (7,0.2075) (8,0.2515) (9,0.2931) (10,0.3326) (11,0.3700) (12,0.4054) (13,0.4388) (14,0.4705) (15,0.5004) (16,0.5287) (17,0.5555) (18,0.5807) (19,0.6046) (20,0.6271) (21,0.6484) (22,0.6684) (23,0.6874) (24,0.7052) (25,0.7221) (26,0.7380) (27,0.7530) (28,0.7672) (29,0.7805) (30,0.7931) (31,0.8050) (32,0.8162) (33,0.8267) (34,0.8367) (35,0.8461) (36,0.8549) (37,0.8633) (38,0.8711) (39,0.8785) (40,0.8855) (41,0.8921) (42,0.8983) (43,0.9042) (44,0.9097) (45,0.9149) (46,0.9198) (47,0.9244) (48,0.9288) (49,0.9329) (50,0.9368) (51,0.9404) (52,0.9438) (53,0.9471) (54,0.9501) (55,0.9530) (56,0.9557) (57,0.9583) (58,0.9607) (59,0.9629) (60,0.9651) (61,0.9671) (62,0.9690) (64,0.9725) (65,0.9740) (66,0.9755) (68,0.9783) (69,0.9795) (71,0.9818) (72,0.9829) (74,0.9848) (75,0.9857) (77,0.9873) (78,0.9880) (80,0.9894) (81,0.9900) (83,0.9911) (85,0.9921) (87,0.9930) (88,0.9934) (90,0.9941) (92,0.9948) (94,0.9954) (96,0.9959) (98,0.9963) (100,0.9968) (102,0.9971) (104,0.9974) (106,0.9977) (108,0.9980) (111,0.9983) (113,0.9985) (115,0.9987) (118,0.9989) (120,0.9990) (122,0.9991) (125,0.9993) (128,0.9994) (130,0.9995) (133,0.9995) (136,0.9996) (138,0.9997) (141,0.9997) (144,0.9998) (147,0.9998) (150,0.9998) (153,0.9999) (156,0.9999) (160,0.9999) (163,0.9999) (166,0.9999) (170,0.9999) (173,1.0000) (177,1.0000) (180,1.0000) (184,1.0000) (188,1.0000) (192,1.0000) (196,1.0000) (200,1.0000) (204,1.0000) (208,1.0000) (212,1.0000) (217,1.0000) (221,1.0000) (226,1.0000) (230,1.0000) (235,1.0000) (240,1.0000) (245,1.0000) (250,1.0000) (255,1.0000) (260,1.0000) (265,1.0000) (271,1.0000) (277,1.0000) (282,1.0000) (288,1.0000) (294,1.0000) (300,1.0000)};
\addlegendentry{$4$}
\addplot[c3, densely dashdotdotted] coordinates {(2,-0.1745) (3,-0.1435) (4,-0.1276) (5,-0.1145) (6,-0.1020) (7,-0.0897) (8,-0.0773) (9,-0.0648) (10,-0.0522) (11,-0.0396) (12,-0.0269) (13,-0.0142) (14,-0.0015) (15,0.0111) (16,0.0237) (17,0.0363) (18,0.0488) (19,0.0612) (20,0.0736) (21,0.0858) (22,0.0980) (23,0.1101) (24,0.1220) (25,0.1339) (26,0.1456) (27,0.1572) (28,0.1687) (29,0.1801) (30,0.1913) (31,0.2025) (32,0.2135) (33,0.2244) (34,0.2351) (35,0.2458) (36,0.2563) (37,0.2666) (38,0.2769) (39,0.2870) (40,0.2970) (41,0.3069) (42,0.3167) (43,0.3263) (44,0.3358) (45,0.3452) (46,0.3545) (47,0.3636) (48,0.3727) (49,0.3816) (50,0.3904) (51,0.3990) (52,0.4076) (53,0.4161) (54,0.4244) (55,0.4326) (56,0.4407) (57,0.4488) (58,0.4566) (59,0.4644) (60,0.4721) (61,0.4797) (62,0.4872) (64,0.5018) (65,0.5090) (66,0.5161) (68,0.5299) (69,0.5367) (71,0.5500) (72,0.5565) (74,0.5692) (75,0.5755) (77,0.5877) (78,0.5936) (80,0.6053) (81,0.6111) (83,0.6223) (85,0.6332) (87,0.6438) (88,0.6489) (90,0.6591) (92,0.6689) (94,0.6785) (96,0.6878) (98,0.6969) (100,0.7056) (102,0.7142) (104,0.7225) (106,0.7305) (108,0.7383) (111,0.7496) (113,0.7569) (115,0.7640) (118,0.7742) (120,0.7808) (122,0.7871) (125,0.7964) (128,0.8052) (130,0.8109) (133,0.8191) (136,0.8269) (138,0.8320) (141,0.8392) (144,0.8462) (147,0.8529) (150,0.8593) (153,0.8654) (156,0.8713) (160,0.8787) (163,0.8839) (166,0.8890) (170,0.8954) (173,0.8999) (177,0.9057) (180,0.9098) (184,0.9150) (188,0.9199) (192,0.9245) (196,0.9289) (200,0.9330) (204,0.9368) (208,0.9405) (212,0.9439) (217,0.9479) (221,0.9509) (226,0.9544) (230,0.9571) (235,0.9601) (240,0.9630) (245,0.9656) (250,0.9681) (255,0.9704) (260,0.9725) (265,0.9745) (271,0.9767) (277,0.9786) (282,0.9802) (288,0.9819) (294,0.9834) (300,0.9848)};
\addlegendentry{$8$}
\end{axis}
\ifdefined\figleftdec\else\newdimen\figleftdec\newdimen\figtmpdim\fi
\pgfextractx{\figleftdec}{\pgfpointanchor{current axis}{west}}
\pgfextractx{\figtmpdim}{\pgfpointanchor{current bounding box}{west}}
\advance\figleftdec by -\figtmpdim
\path (current axis.east) ++(\figleftdec,0);
\end{tikzpicture}
		\end{minipage}
		\caption{Suboptimality index $\alpha_N$ for Ex.~\ref{ex:rates}; dashed: $\alpha_N\in\{0,0.5,0.8\}$.}
		\label{fig:alphaN}
	\end{figure}

	Since Corollary~\ref{cor:discrete} yields exponential cost controllability in general only on balls $\overline{D(M)}_{\bar\rho}$, which are not necessarily invariant under the MPC closed loop, \citet[Thm.~5.2]{Gruene09} does not apply directly, and we only expect asymptotic stability on sublevel sets of $V_N$ for sufficiently large horizons $N$. These horizons can be computed explicitly. Precisely, for $\beta(r,n)=C\sigma^nr$ as in Corollary~\ref{cor:discrete}, we have $\sum_{n=0}^{k-1}\beta(r,n)=\gamma_kr$ with $\gamma_k\coloneqq C\frac{1-\sigma^k}{1-\sigma}$, and the linear program of \citet{Gruene09}, cf.\ also \citet{GruenePannekSeehafer10}, yields the suboptimality index
	\[
	\alpha_N=1-\frac{(\gamma_N-1)\prod_{k=2}^N(\gamma_k-1)}{\prod_{k=2}^N\gamma_k-\prod_{k=2}^N(\gamma_k-1)} .
	\]
	If $\alpha_N>0$, MPC with horizon $N$ is asymptotically stable and its closed-loop cost is at most $V_\infty/\alpha_N$, i.e., at most twice or $1.25$ times the optimal infinite-horizon cost for $\alpha_N\ge0.5$ or $\alpha_N\ge0.8$, respectively. Figure~\ref{fig:alphaN} shows $\alpha_N$ for Example~\ref{ex:rates} with $\lambda=0.1$ and $T=0.5$, that is, $C=1+\frac\lambda2$ and $\sigma=e^{-\min\{1,r_F^2/\bar\rho^2\}}$.

	In the next subsection, we exploit Lemma~\ref{lem:diss} to obtain stability already for $N=2$.

	\subsection{Instantaneous MPC: Stability with horizon two}\label{ssec:N2}
	We now show that MPC with the minimal horizon $N=2$ stabilizes the system on every ball. By Lemma~\ref{lem:diss}, the
	control cost along the feedback \eqref{eq:fb} is bounded by
	$\frac\lambda2\ell^*(x_0)$, which is smaller than $\ell^*(x_0)$ only for
	$\lambda<2$. For larger $\lambda$, we scale the gain,
	\begin{equation}\label{eq:fb-gain}
		\kappa_\eps(x)\coloneqq\PF\big(u^\star-\eps B^*(x-x^\star)\big),\quad\eps\in(0,1],
	\end{equation}
	which trades a smaller control cost for slower decay.
	\begin{lem}\label{lem:gain}
		Let Assumptions~\ref{ass:S} and~\ref{ass:D} hold and let $\eps\in(0,1]$. Then
		$M-B\kappa_\eps$ is maximal monotone with $(M-B\kappa_\eps)(x^\star)=0$. Let
		$(\Scl^\eps(t))_{t\ge0}$ denote the generated semigroup, so that $\Scl^1=\Scl$.
		For $x_0\in\overline{D(M)}$, let $x(t)\coloneqq\Scl^\eps(t)x_0$ and
		$u^\eps_{x_0}(t)\coloneqq\kappa_\eps(x(t))$. Then $u^\eps_{x_0}\in\Uad$
		reproduces the closed loop, and the following holds.
		\begin{enumerate}[label=(\roman*)]
			\item For all $0\le s\le t$,
			\[
			\tfrac12\norm{x(t)-x^\star}^2+\tfrac1\eps\int_s^t\norm{u^\eps_{x_0}-u^\star}_\U^2\,\mathrm d\tau\le\tfrac12\norm{x(s)-x^\star}^2 .
			\]
			\item For all $t\ge0$, $\norm{x(t)-x^\star}^2\le e^{-2\mu^\eps(\rho)t}\rho^2$ with
			$\rho\coloneqq\norm{x_0-x^\star}$ and, with $\mu$ given by \eqref{eq:mudef},
			\begin{equation}\label{eq:mu-gain}
				\mu^\eps(\rho)\coloneqq(1-\eps)\,\zet_M+\eps\,\mu(\rho).
			\end{equation}
		\end{enumerate}
	\end{lem}

	\begin{pf}
		We rescale the input. Let $\tilde B\coloneqq\sqrt\eps B$, $\tilde F\coloneqq\eps^{-1/2}F$ and
		$\tilde u^\star\coloneqq\eps^{-1/2}u^\star\in\tilde F$. Then $(M,\tilde B)$ is a
		maximal monotone control system with controlled equilibrium
		$(x^\star,\tilde u^\star)$. Since $P_{\tilde F}(v)=\eps^{-1/2}\PF(\sqrt\eps v)$,
		its feedback \eqref{eq:fb} is $\tilde\kappa=\eps^{-1/2}\kappa_\eps$, and
		$\tilde B\tilde\kappa=B\kappa_\eps$, so it generates the closed loop of
		\eqref{eq:fb-gain}. It remains to verify Assumptions~\ref{ass:S} and~\ref{ass:D}
		for the rescaled system. By Lemma~\ref{lem:proj-gain} in
		Appendix~\ref{app:proj-gain}, $\norm{\PF(u^\star-\eps\eta)-u^\star}\ge\eps\norm{\PF(u^\star-\eta)-u^\star}$
		for all $\eta\in\U$. With $\eta=B^*(x-x^\star)$, $\norm{x-x^\star}\le\rho$ and
		Assumption~\ref{ass:S}, this gives
		\begin{align*}
			\norm{P_{\tilde F}(\tilde u^\star-\tilde B^*(x-x^\star))-\tilde u^\star}_\U \ge\eps^{-1/2}\eps\,\theta(\rho)\norm{\eta}_\U 
			=\theta(\rho)\norm{\tilde B^*(x-x^\star)}_\U .
		\end{align*}
		Thus, Assumption~\ref{ass:S} holds for $(M,\tilde B)$ with the same $\theta$. Further, the left-hand side of \eqref{eq:D} for $\tilde B$ is
		\begin{align*}
			(1-\eps)\ip{M(x)-M(x^\star)}{x-x^\star}+\eps\big(\ip{M(x)-M(x^\star)}{x-x^\star}+\norm{B^*(x-x^\star)}_\U^2\big),
		\end{align*}
		which is at least $\zet^\eps(\rho)\norm{x-x^\star}^2$ with
		$\zet^\eps\coloneqq(1-\eps)\zet_M+\eps\zet\ge\zet_M$ by \eqref{eq:DM} and
		\eqref{eq:D}. Thus, Assumption~\ref{ass:D} holds with the same $\zet_M$. The
		corresponding rate \eqref{eq:mudef} is
		\[
			(1-\theta^2)\zet_M+\theta^2\zet^\eps =(1-\eps\theta^2)\zet_M+\eps\theta^2\zet =(1-\eps)\zet_M+\eps\big((1-\theta^2)\zet_M+\theta^2\zet\big)=\mu^\eps(\rho).
		\]
		Finally, \citet[Lems.~3.2, 3.3]{PGS26} and Lemma~\ref{lem:adm} for the rescaled
		system give the first claims, Lemma~\ref{lem:diss} together with
		$\tilde\kappa-\tilde u^\star=\eps^{-1/2}(\kappa_\eps-u^\star)$ gives (i), and
		Theorem~\ref{thm:coercive} gives (ii).
	\end{pf}

	We now provide a bound on the value function. Let $J_k(x_0,u)\coloneqq\sum_{n=0}^{k-1}l(x_u(n),u(n))$ for $k\in\mathbb{N}\cup\{\infty\}$ and $u\in\Useq$, so that $V_k=\inf_uJ_k$ for $k\in\mathbb{N}$ and $V_\infty\coloneqq\inf_uJ_\infty$.
	\begin{lem}\label{lem:accum}
		Let Assumptions~\ref{ass:S} and~\ref{ass:D} hold, let $\lambda>0$, $\bar\rho>0$ and
		$\eps\in(0,1]$ with $\lambda\eps<2$. Set $\sigma^\eps\coloneqq e^{-2\mu^\eps(\bar\rho)T}$
		with $\mu^\eps$ from \eqref{eq:mu-gain} and
		\[
		\hat\gamma_k^\eps\coloneqq\tfrac{1-(\sigma^\eps)^{k-1}}{1-\sigma^\eps}+\tfrac{\lambda\eps}2
		+\big(1-\tfrac{\lambda\eps}2\big)(\sigma^\eps)^{k-1},\qquad k\in\mathbb{N}.
		\]
		Then, for all $x_0\in\overline{D(M)}_{\bar\rho}$ and all $k\in\mathbb{N}$,
		\[
		V_k(x_0)\ \le\ \hat\gamma_k^\eps\,l^*(x_0)\ \le\ \Big(\tfrac{1}{1-\sigma^\eps}+\tfrac{\lambda\eps}2\Big)\,l^*(x_0).
		\]
	\end{lem}
	\begin{pf}
		By Lemma~\ref{lem:gain}(ii) and since $\mu^\eps$ is non-increasing, we have for all $x\in\overline{D(M)}_{\bar\rho}$,
		\begin{equation}\label{eq:contr}
			\norm{\Scl^\eps(T)x-x^\star}^2\le\sigma^\eps\,\norm{x-x^\star}^2 .
		\end{equation}
		For $k=1$, the claim holds with equality, since $V_1=l^*$ by definition. Let $k\ge2$ and $x_0\in\overline{D(M)}_{\bar\rho}$,
		and let
		\begin{equation*}
			\tilde u(n)\coloneqq \begin{cases}
				u^\eps_{x_0}(nT+\cdot)|_{[0,T)}, & n=0,\dots,k-2 ; \\
				u^\star, & n\ge k-1,
			\end{cases}
		\end{equation*}
		with $u^\eps_{x_0}$ from Lemma~\ref{lem:gain}, which is
		admissible by Lemma~\ref{lem:gain}. We abbreviate
		$a_n\coloneqq\norm{x_{\tilde u}(n)-x^\star}^2$. For $n\le k-1$, we
		have $x_{\tilde u}(n)=\Scl^\eps(nT)x_0$, such that Lemma~\ref{lem:gain}(i) keeps the state in
		$\overline{D(M)}_{\bar\rho}$ and \eqref{eq:contr}, applied successively, yields
		$a_n\le(\sigma^\eps)^na_0$. Moreover, Lemma~\ref{lem:gain}(i) with $s=0$ and $t=(k-1)T$ yields
		$$\sum_{n=0}^{k-2}\int_{nT}^{(n+1)T}\norm{u^\eps_{x_0}(\tau)-u^\star}^2\,\mathrm d\tau\le\frac{\eps}{2}(a_0-a_{k-1}).$$ Since the last stage is evaluated at $u^\star$,
		that is, $l(x_{\tilde u}(k-1),u^\star)=a_{k-1}$, we obtain
		\[
		\begin{aligned}
			V_k(x_0)&\le J_k(x_0,\tilde u)\\
			&=\sum_{n=0}^{k-1}a_n+\lambda\sum_{n=0}^{k-2}\int_{nT}^{(n+1)T}\norm{u^\eps_{x_0}(\tau)-u^\star}^2\,\mathrm d\tau\\
			&\le\sum_{n=0}^{k-2}a_n+\big(1-\tfrac{\lambda\eps}2\big)a_{k-1}+\tfrac{\lambda\eps}2a_0 \leq \hat\gamma_k^\eps\,a_0 ,
		\end{aligned}
		\]
		where the last inequality uses $a_n\le(\sigma^\eps)^na_0$ and $\lambda\eps<2$, and $a_0=l^*(x_0)$ by \eqref{eq:lstar}, which shows the claim.
	\end{pf}

	Since $\eps$ only enters the analysis, we choose $\eps=\min\{1,1/\lambda\}$, which, for $\zet_M=0$ and $T\mu(\bar\rho)\ll1$, approximately minimizes the contraction factor $\hat\gamma_2^\eps-1=1-(1-\frac{\lambda\eps}{2})(1-\sigma^\eps)$ obtained below, as $(1-\frac{\lambda\eps}{2})(1-\sigma^\eps)\approx2T\mu(\bar\rho)\,\eps(1-\frac{\lambda\eps}{2})$.

	\begin{thm}[Stabilization with horizon two]\label{thm:N2}
		Let Assumptions~\ref{ass:S} and~\ref{ass:D} hold, let $\lambda>0$ be the control penalization in the running cost \eqref{eq:l}, let $T>0$ be the sampling time and consider
		the MPC scheme \eqref{eq:ocp} with $N=2$. For $\bar\rho>0$, let
		$\sigma_\lambda\coloneqq\exp\Big(-2T\Big(\zet_M+\tfrac{\mu(\bar\rho)-\zet_M}{\max\{1,\lambda\}}\Big)\Big)$
		with $\mu$ from \eqref{eq:mudef}, that is, $\sigma_\lambda=\sigma^\eps$ from Lemma~\ref{lem:accum} with $\eps=\min\{1,1/\lambda\}$, and
		\[
		\begin{aligned}
			\hat\gamma_2&=1+\tfrac{\min\{1,\lambda\}}2+\Big(1-\tfrac{\min\{1,\lambda\}}2\Big)\sigma_\lambda,\\
			\alpha_2&\coloneqq1-(\hat\gamma_2-1)^2\in(0,1] .
		\end{aligned}
		\]
		Then, for every $x\in\overline{D(M)}_{\bar\rho}$, the MPC successor
		$x^+\coloneqq f(x,\mu_2(x))$ satisfies the following.
		\begin{enumerate}[label=(\roman*)]
			\item $\norm{x^+-x^\star}^2\le(\hat\gamma_2-1)\norm{x-x^\star}^2$, where
			$\hat\gamma_2-1=1-(1-\frac{\min\{1,\lambda\}}{2})(1-\sigma_\lambda)<1$. In particular, $\overline{D(M)}_{\bar\rho}$ is
			forward invariant under the MPC closed loop, and
			$$\norm{x_{\mu_2}(n)-x^\star}^2\le(\hat\gamma_2-1)^n\norm{x_0-x^\star}^2$$ for all
			$x_0\in\overline{D(M)}_{\bar\rho}$ and $n\in\mathbb{N}_0$, that is, as $\bar\rho>0$ is arbitrary, $x^\star$ is globally
			asymptotically stable with semiglobal exponential rate for the MPC closed loop.
			\item $V_2(x^+)\le V_2(x)-\alpha_2\,l(x,\mu_2(x))$ and, with the closed-loop cost
			$J^{\mathrm{cl}}_\infty(x,\mu_2)\coloneqq J_\infty(x,u)$ with $u(n)=\mu_2(x_{\mu_2}(n))$ and $x_{\mu_2}(0)=x$, we have the suboptimality estimate
			$$J^{\mathrm{cl}}_\infty(x,\mu_2)\le V_2(x)/\alpha_2\le V_\infty(x)/\alpha_2.$$
		\end{enumerate}
	\end{thm}

	\begin{pf}
		Let $x\in\overline{D(M)}_{\bar\rho}$, let $(\hat u(0),\hat u(1))$ be an optimal control sequence for
		$V_2(x)$, such that $\mu_2(x)=\hat u(0)$ and $x^+=x_{\hat u}(1)$, and set
		$ l_0\coloneqq l(x,\hat u(0))$ and $ l_1\coloneqq l(x^+,\hat u(1))$. By
		Lemma~\ref{lem:accum} with $k=2$ and $\eps=\min\{1,1/\lambda\}$, we have
		$ l_0+ l_1=V_2(x)\le\hat\gamma_2\,l^*(x)$. Since $l^*(x)\le l_0$ and
		$l^*(x^+)\le l_1$ by definition of $l^*$, this yields, together with \eqref{eq:lstar},
		$\norm{x^+-x^\star}^2=l^*(x^+)\le l_1\le(\hat\gamma_2-1)\,l^*(x)$, which shows the first
		estimate in (i), and the remaining assertions of (i) follow by iteration.
		Hence,
		$x^+\in\overline{D(M)}_{\bar\rho}$ and Lemma~\ref{lem:accum} at $x^+$ gives
		$V_2(x^+)\le\hat\gamma_2\,l^*(x^+)\le\hat\gamma_2 l_1$.
		Then, the performance estimate follows from the relaxed Lyapunov inequality
		\[
			V_2(x)-V_2(x^+) \ge l_0-(\hat\gamma_2-1) l_1 \ge l_0-(\hat\gamma_2-1)^2\,l^*(x)\ge\alpha_2 l_0 ,
		\]
		due to~\citet[Thm.~4.11]{GruenePannek}.
	\end{pf}

	For $\lambda>1$, the guaranteed contraction deteriorates as $\lambda$ grows. Theorem~\ref{thm:N2} covers $u^\star\in\operatorname{int}F$, e.g., $L^2$-norm bounds, but not the pointwise constraints of Example~\ref{ex:heat}, for which Section~\ref{ssec:numerics} nevertheless indicates asymptotic stability. This is plausible, as $L^\infty(\Omega)$-balls are invariant under the closed loop of \eqref{eq:fb} by Lemma~\ref{lem:heatLinf}, such that, for $\lambda<2$, the first estimate of Theorem~\ref{thm:N2}(i) holds on them with $\min\{1,\lambda\}$ replaced by $\lambda$ and $\sigma_\lambda$ by $e^{-2\mu T}$, $\mu$ from Example~\ref{ex:heat}. Its iteration would require their invariance under the MPC closed loop.

	\begin{cor}[Horizon two under exponential detectability]\label{cor:N2detect}
		Let the hypotheses of Theorem~\ref{thm:detect} hold, let $\lambda>0$ be the control penalization in the running cost \eqref{eq:l}, let $\eps\in(0,1]$ with $\lambda\eps<2$ and $\bar\rho>0$, let $K^\eps$ and $\nu^\eps$ be defined as $K$ and $\nu$ in Theorem~\ref{thm:detect} with $g$ replaced by $g^\eps(\rho)\coloneqq\sqrt\eps\norm{B}+\norm{L}/(\sqrt\eps\,\theta(\rho))$, and let the sampling time satisfy $T>4K^\eps(\bar\rho)^2$. Then the assertions of Lemma~\ref{lem:accum} and of Theorem~\ref{thm:N2} with $\min\{1,\lambda\}$ replaced by $\lambda\eps$ hold on $\overline{D(M)}_{\bar\rho}$ with $\sigma^\eps$ and $\sigma_\lambda$, respectively, replaced by $\sigma\coloneqq4e^{-2\nu^\eps(\bar\rho)T}\in(0,1)$, except that, as $T$ depends on $\bar\rho$, the global asymptotic stability in Theorem~\ref{thm:N2}(i) reduces to exponential stability on $\overline{D(M)}_{\bar\rho}$.
	\end{cor}

	\begin{pf}
		Let $\tilde B=\sqrt\eps B$ be as in the proof of Lemma~\ref{lem:gain}, which shows that $(M,\tilde B)$ satisfies Assumption~\ref{ass:S} with the same $\theta$ and has the closed loop $\Scl^\eps$. Since $\tilde L\coloneqq L/\sqrt\eps$ satisfies $\tilde L\tilde B^*=LB^*$, $\norm{\tilde L}=\norm{L}/\sqrt\eps$ and $\norm{\tilde B}=\sqrt\eps\norm{B}$, Theorem~\ref{thm:detect} applies to $(M,\tilde B)$ with $g^\eps$ in place of $g$. As $\nu^\eps$ is non-increasing, \eqref{eq:detect-decay} for $\Scl^\eps$ yields \eqref{eq:contr} with $\sigma$, and $\sigma<1$ if and only if $T>\ln2/\nu^\eps(\bar\rho)=4K^\eps(\bar\rho)^2$. The claim follows, as the proofs of Lemma~\ref{lem:accum} and Theorem~\ref{thm:N2} use Assumptions~\ref{ass:S} and~\ref{ass:D} only via \eqref{eq:contr}, while the remaining claims of Lemma~\ref{lem:gain} hold without them.
	\end{pf}

	\subsection{Numerical illustration}\label{ssec:numerics}
\begin{figure}[t]
		\centering
		\begin{minipage}{\figwidth}\centering
			\input{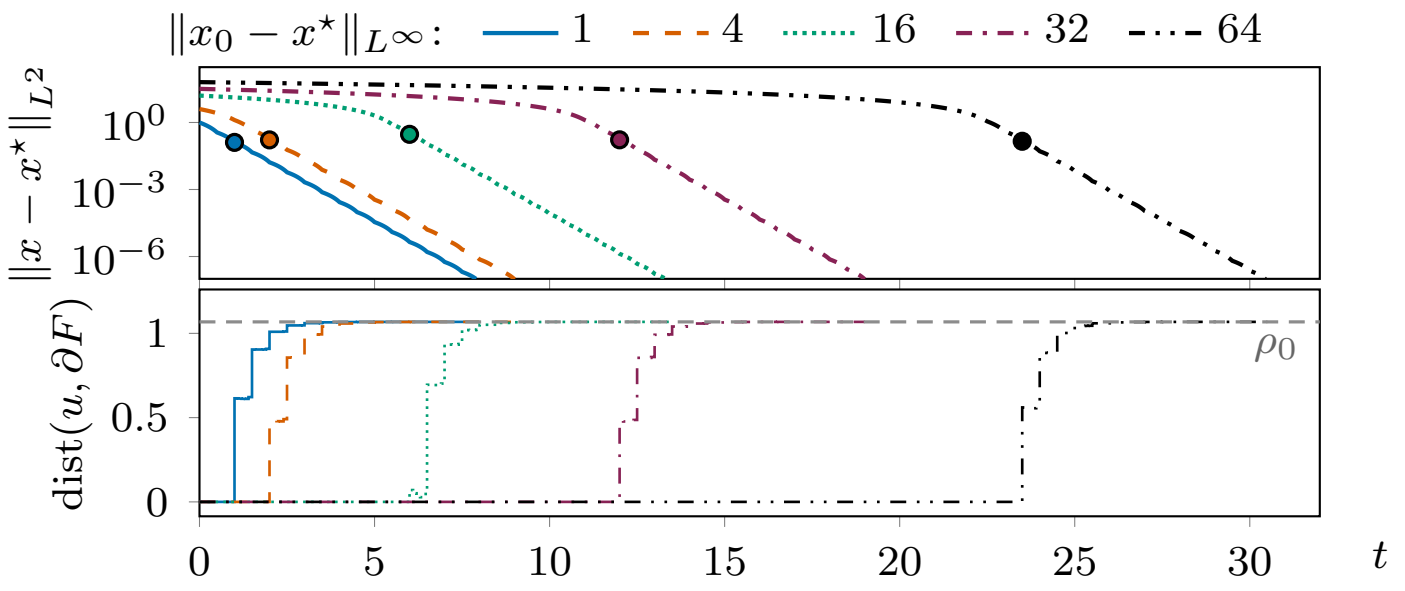}
		\end{minipage}
		\caption{MPC closed loop for Ex.~\ref{ex:heat} with
			$f(s)=s^3/(1+|s|^3)$, $\Omega_c=(0.2,0.8)^2$, $u^\star=4\sin(2\pi\xi)\sin(2\pi\eta)$ on $\Omega_c$,
			$(T_{\min},T_{\max})=(-5,7)$,
			$x_0-x^\star\equiv R_\infty\in\{1,4,16,32,64\}$. }
		\label{fig:mpc-heat}
	\end{figure}

\begin{figure}[t]
		\centering
		\begin{minipage}{\figwidth}\centering
			\input{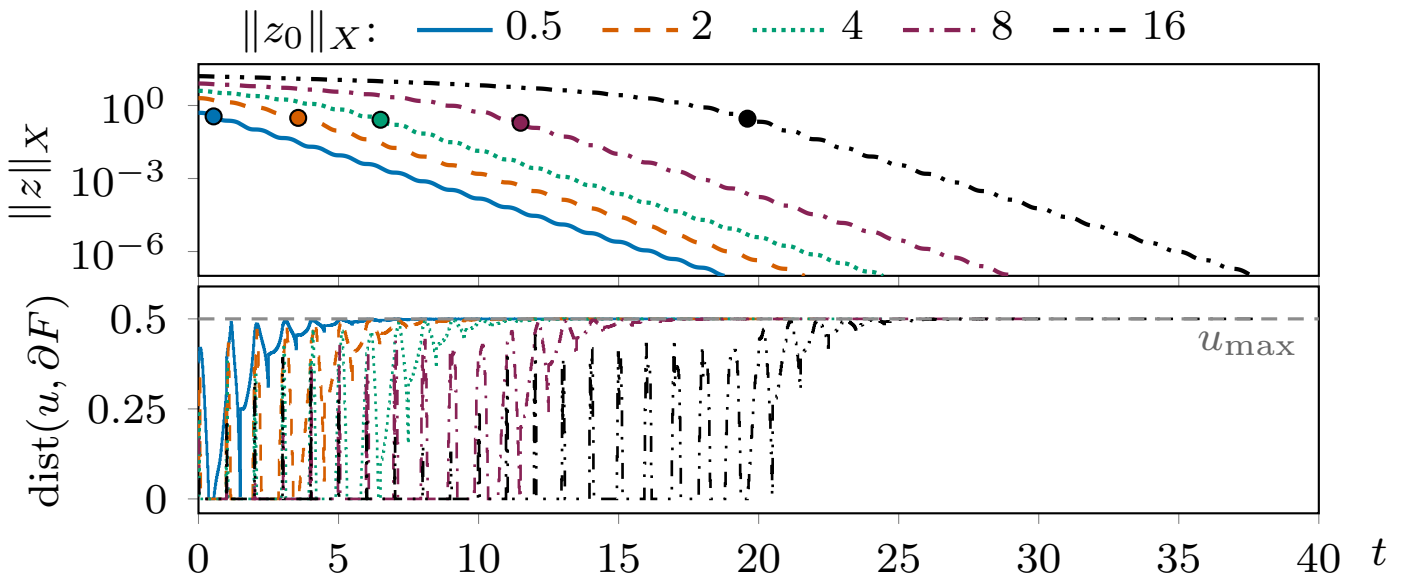}
		\end{minipage}
		\caption{MPC closed loop for Ex.~\ref{ex:wave} on
			$\Omega=(0,1)$, with $\psi=\arctan$,
			$\Omega_c=(0.1,0.3)$, $u_{\max}=0.5$ and initial displacements proportional to
			$\sin(\pi\xi)$ with $\norm{z_0}_\X\in\{0.5,2,4,8,16\}$. }
		\label{fig:mpc-wave}
	\end{figure}

	We illustrate MPC with $N=2$, $T=0.5$ and $\lambda=0.1$ for Examples~\ref{ex:heat} and~\ref{ex:wave}; the figures for Example~\ref{ex:rates} look similar and are omitted. The controls are piecewise constant in time, and \eqref{eq:ocp} is solved by gradient-based methods. Figures~\ref{fig:mpc-heat} and~\ref{fig:mpc-wave} show the closed loop for increasing initial deviations and the distance of the control to $\partial F$. For large deviations, the control stays at the boundary for a time growing with the deviation, during which the state decays slowly, and afterwards all trajectories decay with the same exponential rate. This reflects the semiglobal exponential rates in Example~\ref{ex:heat} and Theorem~\ref{thm:detect}.

	\section{Conclusion}\label{sec:conclusion}

	Using a saturated feedback as admissible control, we have shown that maximal monotone control systems under input constraints are globally cost controllable with semiglobal rate. A dissipation inequality reduces cost controllability to a decay estimate for the closed-loop state, for which we established explicit semiglobal exponential rates under coercivity of $M+BB^*$ or exponential detectability of the linear part. In the sampled-data setting under coercivity, MPC without terminal conditions and horizon $N=2$ renders the system globally asymptotically stable with semiglobal exponential rate.

	\section*{Declaration of generative AI use}
	During the preparation of this work, the authors used Claude (Anthropic) to improve language, draft and check arguments, and implement the numerical illustrations.
	The authors take full responsibility for all content.

	\bibliographystyle{plainnat}
	\bibliography{cost-controllability}

	\appendix

	\section{Auxiliary results for the heat equation}\label{app:aux}

	\begin{prop}\label{prop:AN+F}
		Let $\Omega\coloneqq (0,1)^{2}$ and $H\coloneqq L^{2}(\Omega)$. Let
		$f\colon\mathbb{R}\to\mathbb{R}$ be continuous and nondecreasing with $f(0)=0$, and put
		\[
		\mathcal Ny\coloneqq f\circ y,\qquad  D(\mathcal N)\coloneqq \{y\in H: f\circ y\in H\},
		\]
		\[
		\begin{aligned}
			A_{N}y&\coloneqq -\Delta y,\\
			D(A_{N})&\coloneqq \Bigl\{y\in H^{2}(\Omega):\tfrac{\partial y}{\partial\nu}=0 \text{ on }\partial\Omega\Bigr\}.
		\end{aligned}
		\]
		Then $D(A_{N})\subseteq D(\mathcal N)$ and $A_{N}+\mathcal N$ with domain $D(A_{N})$ is
		maximal monotone in $H$.
	\end{prop}

	\begin{pf}
		Let $\phi(s)\coloneqq \int_{0}^{s}f(r)\,dr$, a convex $C^{1}$ function with $\phi\ge\phi(0)=0$,
		and let
		\[
		\begin{aligned}
			\Phi(y)&\coloneqq \int_{\Omega}\phi(y(x))\,dx ,\\
			D(\Phi)&=\{y\in H:\phi(y)\in L^{1}(\Omega)\}.
		\end{aligned}
		\]
		By \citet[Prop.~2.7]{Barbu}, $\Phi$ is convex, l.s.c.\ and proper on $H$
		(properness from $\Phi(0)=0$), and $\partial\Phi=\mathcal N$.

		To prove the maximal monotonicity we verify the assumptions of
		\citet[Thm.~2.11]{Barbu} with $h=0$, $C=0$. To this end, let $y\in D(\Phi)$,
		$\lambda>0$ and $z\coloneqq (I+\lambda A_{N})^{-1}y \in D(A_N)\subset H^2(\Omega)$, i.e.\ $z-\lambda\Delta z=y$ in
		$\Omega$ and $\partial z/\partial\nu=0$ on $\partial\Omega$; we have to show
		\begin{equation}\label{eq:inv}
			\Phi(z)\le\Phi(y).
		\end{equation}
		Let $\mu>0$, let $f_{\mu} \coloneqq \tfrac{1}{\mu}(I-(I+\mu f)^{-1})$ be the \textit{Yosida approximation} of $f$ and let ${}^{\mu\!}\phi : \mathbb{R} \to \mathbb{R}$ given by
		\begin{equation*}
			{}^{\mu\!}\phi(s) \coloneqq \inf_{r \in \mathbb{R}} (\phi(r) + \tfrac{1}{2\mu} |s-r|^2)
		\end{equation*}
		be the \textit{Moreau envelope} of $\phi$. By \citet[Prop. 12.5]{BausComb2011} ${}^{\mu\!}\phi$ is convex and differentiable such that \citet[Prop. 12.30]{BausComb2011} yields
		\begin{equation*}
			\begin{aligned}
				{}^{\mu\!}\phi'&= \tfrac{1}{\mu}(I-\operatorname{prox}_{\mu \phi}) = \tfrac{1}{\mu}(I-(I+\mu \partial \phi)^{-1})\\
				&=\tfrac{1}{\mu}(I-(I+\mu f)^{-1})=  f_{\mu}.
			\end{aligned}
		\end{equation*}
		Moreover, $0\le{}^{\mu\!}\phi\le\phi$ with ${}^{\mu\!}\phi\uparrow\phi$ pointwise as $\mu\downarrow0$, \citep[Prop. 12.33]{BausComb2011}. Put ${}^{\mu}\Phi(v)\coloneqq \int_{\Omega}{}^{\mu\!}\phi(v)\,dx$. As $f_{\mu}$ is Lipschitz and nondecreasing with $f_{\mu}(0)=0$, \citet[Thm.~2.1.11]{Ziemer} implies $f_{\mu}(z)\in H^{1}(\Omega)$ and it holds that
		\begin{equation*}
			\nabla f_{\mu}(z)=f_{\mu}'(z)\nabla z
		\end{equation*}
		a.e. Moreover, $f_{\mu}'\ge0$ a.e. Multiplying $z-\lambda\Delta z=y$ by $f_{\mu}(z)$ and
		using Green's formula together with $\partial z/\partial\nu=0$,
		\begin{equation}\label{eq:test}
			\begin{aligned}
				\int_{\Omega}f_{\mu}(z)(z-y)\,dx &=\lambda\int_{\Omega}\Delta z\,f_{\mu}(z)\,dx\\
				&=-\lambda\int_{\Omega}f_{\mu}'(z)\,|\nabla z|^{2}\,dx\ \le\ 0 .
			\end{aligned}
		\end{equation}
		By convexity of ${}^{\mu\!}\phi$ we have
		\begin{equation*}
			{}^{\mu\!}\phi(y)-{}^{\mu\!}\phi(z)\ge  {}^{\mu\!}\phi'(z)(y-z)= f_{\mu}(z)(y-z)
		\end{equation*}
		a.e.\ in $\Omega$. A combination of the definition of $ {}^{\mu}\Phi$ and \eqref{eq:test} yields
		\[
		0 \le  \int_{\Omega}f_{\mu}(z)(y-z)\,dx \le \int_{\Omega} {}^{\mu\!}\phi(y)-{}^{\mu\!}\phi(z) \,dx
		\]
		and thus ${}^{\mu}\Phi(z)\ \le\ {}^{\mu}\Phi(y)\ \le\ \Phi(y)$, where the second inequality uses ${}^{\mu\!}\phi\le\phi$. Finally $0\le{}^{\mu\!}\phi(z)\uparrow\phi(z)$ a.e.\ in $\Omega$ as $\mu\downarrow0$, so
		${}^{\mu}\Phi(z)\to\Phi(z)$ by monotone convergence, and \eqref{eq:inv} follows.

		Clearly, $A_{N}$ is maximal monotone in $H$ and $0\in D(A_{N})\cap D(\partial\Phi)$. By \citet[Thm.~2.11]{Barbu}, $A_{N}+\partial\Phi=A_{N}+\mathcal N$ is maximal monotone with $D(A_{N}+\mathcal N)=D(A_{N})\cap D(\mathcal N)$. Since $\Omega\subset\mathbb R^2$, the Sobolev embedding gives $D(A_{N})\subseteq H^{2}(\Omega)\subseteq L^{\infty}(\Omega)$; as $|\Omega|=1$ and $f$ is bounded on compact subsets, $L^{\infty}(\Omega)\subseteq D(\mathcal N)$, so $D(A_{N}+\mathcal N)=D(A_{N})$.
	\end{pf}

	\begin{lem}[$L^\infty$-invariance for the heat equation]\label{lem:heatLinf}
		In the setting of Example~\ref{ex:heat}, for every
		$x_0\in\overline{D(M)}$ with $x_0-x^\star\in L^\infty(\Omega)$,
		\[
		\norm{\Scl(t)x_0-x^\star}_{L^\infty(\Omega)}\le\norm{x_0-x^\star}_{L^\infty(\Omega)}\qquad\text{for all }t\ge0 .
		\]
	\end{lem}

	\begin{pf}
		We show that the sets
		$C_R\coloneqq \{v\in L^2(\Omega):\norm{v}_{L^\infty(\Omega)}\le R\}$, $R>0$, are invariant under
		the shifted closed loop. To this end, we write $w\coloneqq x-x^\star$ and, using the equilibrium
		relation $A_Nx^\star+f(x^\star)=Bu^\star$ and the pointwise action of $\PF$, the shifted closed-loop operator as $W(w)\coloneqq \Mcl(x^\star+w)=A_Nw+b(\cdot,w)$ with
		\[
			b(\xi,s) \coloneqq f\big(x^\star(\xi)+s\big)-f\big(x^\star(\xi)\big)
			+\mathbf 1_{\Omega_c}(\xi)\Big(u^\star(\xi)-P_{[T_{\min},T_{\max}]}\big(u^\star(\xi)-s\big)\Big),
		\]
		where $P_{[T_{\min},T_{\max}]}$ denotes the scalar projection onto $[T_{\min},T_{\max}]$.
		For each $\xi$, the function $b(\xi,\cdot)$ is non-decreasing with $b(\xi,0)=0$: the
		first summand due to the monotonicity of $f$, the second since
		$s\mapsto u^\star-P_{[T_{\min},T_{\max}]}(u^\star-s)$ is non-decreasing and vanishes at
		$s=0$ due to $u^\star\in F$. Now let $g \in  C_R$, $h>0$ and define
		$w \coloneqq (I+h W)^{-1}g\in D(A_N)$, that is,
		\begin{equation}\label{eq:resolvent_identity}
			w+h A_Nw+h b(\cdot,w)=g.
		\end{equation}
		Note that $g \in C_R$ means $-R \leq g(\xi) \leq R$ for a.e. $\xi \in \Omega$. Testing \eqref{eq:resolvent_identity} with the truncation $(w-R)^+\in H^1(\Omega)$ and using the form identity
		$\ip{A_Nw}{\varphi}_{L^2(\Omega)}=\int_\Omega\nabla w\cdot\nabla\varphi\,d\xi$ for
		$\varphi\in H^1(\Omega)$, we obtain
		\begin{equation*}
			\int_\Omega (w-g)(w-R)^+ +h\int_\Omega\nabla w\cdot\nabla (w-R)^+ +h\int_\Omega b(\cdot, w)(w-R)^+ =0 .
		\end{equation*}
		The diffusion term $h\int_\Omega\nabla w\cdot\nabla (w-R)^+$ is nonnegative, and on $\{w>R\}$ we have $b(\cdot,w) \geq b(\cdot,0) = 0$. Thus,
		$0\ge\int_\Omega (w-g)(w-R)^+\ge\int_\Omega \big((w-R)^+\big)^2$, where we used $w-g\ge w-R$. Consequently, $w\le R$. Analogously, testing \eqref{eq:resolvent_identity} with $-(w+R)^-$
		yields $w\ge-R$. The resolvent of $W$ thus maps $C_R$ into itself, and the
		exponential formula
		$\Scl(t)x_0-x^\star=\lim_{n\to\infty}\big(I+\tfrac tnW\big)^{-n}(x_0-x^\star)$, cf.\
		\citep{Barbu}, together with the convexity and $L^2$-closedness of $C_R$, \citep[Thm. 4.9]{Brezis2011FunctionalAnalysis}, yields
		the assertion for $R=\norm{x_0-x^\star}_{L^\infty(\Omega)}$.
	\end{pf}

	\section{Proof of Theorem~\ref{thm:detect}}\label{app:detect}

	The proof relies on the following lemma, which turns an algebraic decay estimate into an exponential one.

	\begin{lem}\label{lem:alg2exp}
		Let $\rho>0$, let $K\colon[0,\rho]\to(0,\infty)$ be non-decreasing, and let
		$d\colon[0,\infty)\to[0,\rho]$ be non-increasing with $d(0)\le\rho$ such that, for all $s\ge0$ and $t>0$,
		\[
		d(s+t)\le K(d(s))\,d(s)/\sqrt t .
		\]
		Then, with $\nu(\rho)\coloneqq\ln2/(4K(\rho)^2)$, for all $t\ge0$,
		\[
		d(t)\le 2e^{-\nu(\rho)t}\,d(0).
		\]
	\end{lem}

	\begin{pf}
		Set $t_0\coloneqq 4K(\rho)^2$. We claim $d(kt_0)\le2^{-k}d(0)$ for all $k\in\mathbb{N}_0$ and proceed by induction. The case $k=0$ is trivial. Assume the claim for $k$ and put $\rho_k\coloneqq d(kt_0)\le2^{-k}d(0)\le\rho$. Since $K$ is non-decreasing, $K(\rho_k)\le K(\rho)$, such that $t_k\coloneqq4K(\rho_k)^2\le t_0$, and the assumption with $s=kt_0$ and $t=t_k$ gives $d(kt_0+t_k)\le\rho_k/2$. As $d$ is non-increasing,
		$d((k+1)t_0)\le d(kt_0+t_k)\le\rho_k/2\le2^{-(k+1)}d(0)$, which proves the claim.
		For $t\ge0$, monotonicity of $d$ and $\lfloor\tau\rfloor>\tau-1$ yield
		\begin{align*}
			d(t)  \le d(\lfloor t/t_0\rfloor t_0)  \le2^{-\lfloor t/t_0\rfloor}d(0) \le2\cdot2^{-t/t_0}d(0) =2e^{-\nu(\rho)t}d(0),
		\end{align*}
		since $\ln2/t_0=\nu(\rho)$.
	\end{pf}

	\begin{proof}[Proof of Theorem~\ref{thm:detect}]
		We fix $x_0\in D(M)$, as the general case follows by density as in the proof of
		Theorem~\ref{thm:coercive}, using continuity of $\nu$.  Set $\rho\coloneqq\norm{x_0-x^\star}$ and abbreviate, as in the proof of Theorem~\ref{thm:coercive},
		$x(t)\coloneqq\Scl(t)x_0$, $\delta x(t)\coloneqq x(t)-x^\star$, $\delta y(t)\coloneqq B^*\delta x(t)$ and $\delta u(t)\coloneqq\PF(u^\star-\delta y(t))-u^\star$, as well as
		$\Delta N(t)\coloneqq N(x(t))-N(x^\star)$.
		By \eqref{eq:dissipation-core}, \eqref{eq:key-sector}, the monotonicity of $A$ and
		\eqref{eq:cocoercive}, we have
		$\tfrac{\mathrm d}{\mathrm{d}t}\tfrac12\norm{\delta x}^2\le-c_N\norm{\Delta N}^2-\norm{\delta u}^2$
		for a.e.\ $t\ge0$. Integration yields
		\begin{equation*}
			\tfrac12\norm{\delta x(t)}^2+\int_0^t\big(c_N\norm{\Delta N(\tau)}^2+\norm{\delta u(\tau)}_\U^2\big)\,\mathrm d\tau \le \tfrac12\rho^2 .
		\end{equation*}
		By dropping the first term and letting $t \to \infty$:
		\[
		c_N\norm{\Delta N}_{L^2}^2+\norm{\delta u}_{L^2}^2\le\tfrac12\rho^2,
		\]
		where $L^2\coloneqq L^2(0,\infty)$. Since $$-A\delta x-\Delta N+B\delta u=(-A+LB^*)\delta x+B\delta u-L\delta y-\Delta N,$$ the
		variation-of-parameters formula, cf.\ \citep[Eq. 11.3]{CurtainZwart}, yields with $T(t) \coloneqq T_{-A+LB^*}(t)$:
		\[
			\delta x(t)=T(t)\delta x(0) +\int_0^t T(t-\tau)\big(B\delta u(\tau)-L\delta y(\tau)-\Delta N(\tau)\big)\,\mathrm{d}\tau .
		\]
		Since $\norm{\delta x(\tau)}\le\rho$, Assumption~\ref{ass:S} implies
		$\norm{\delta y(\tau)}\le\theta(\rho)^{-1}\norm{\delta u(\tau)}$ and hence
		$\norm{B\delta u(\tau)-L\delta y(\tau)}\le g(\rho)\norm{\delta u(\tau)}$ for all $\tau\ge0$.
		Thus
		\begin{equation}\label{eq:vop}
				\norm{\delta x(t)} \le c\,e^{-\omega t}\rho+c\int_0^t e^{-\omega(t-\tau)}\big(g(\rho)\norm{\delta u(\tau)}
				+\norm{\Delta N(\tau)}\big)\,\mathrm{d}\tau .
		\end{equation}
		Taking $L^2$-norms in \eqref{eq:vop}, Young's convolution inequality,
		\citep[Thm. 4.15]{Brezis2011FunctionalAnalysis}, $\norm{e^{-\omega\cdot}}_{L^2}=(2\omega)^{-1/2}$,
		$\norm{e^{-\omega\cdot}}_{L^1}=\omega^{-1}$ as well as $\norm{\delta u}_{L^2}\le\rho/\sqrt2$ and
		$\norm{\Delta N}_{L^2}\le\rho/\sqrt{2c_N}$ yield
		\begin{align*}
			\norm{\delta x}_{L^2}&\le\tfrac{c\rho}{\sqrt{2\omega}}+\tfrac{c}{\omega}\big(g(\rho)\norm{\delta u}_{L^2}+\norm{\Delta N}_{L^2}\big)\\
			&\le \tfrac{c\rho}{\sqrt{2\omega}}+\tfrac{c\,(g(\rho)+c_N^{-1/2})\rho}{\sqrt2\,\omega}=K(\rho)\rho .
		\end{align*}
		As $\tau\mapsto\norm{\delta x(\tau)}$ is non-increasing by Lemma~\ref{lem:diss}, we have
		$t\norm{\delta x(t)}^2\le\int_0^t\norm{\delta x(\tau)}^2\,\mathrm{d}\tau\le K(\rho)^2\rho^2$,
		that is,
		\begin{equation}\label{eq:alg}
			\norm{\delta x(t)}\le K(\rho)\rho/\sqrt t\quad\text{for all }t>0 .
		\end{equation}
		Since $\Scl$ is a semigroup, the estimate \eqref{eq:alg} holds for every initial value $\Scl(s)x_0$ with $\rho$ replaced by $\norm{\delta x(s)}\le\rho$. Hence, Lemma~\ref{lem:alg2exp} yields \eqref{eq:detect-decay}.
		Finally, $\theta$ is continuous and non-increasing by
		Assumption~\ref{ass:S}, hence $g$ and $K$ are continuous and non-decreasing and $\nu$ is
		continuous and non-increasing. As in the proof of Theorem~\ref{thm:coercive}, this yields
		$\tilde\beta(r,t)=4e^{-2\nu(\sqrt r)t}r\in\mathcal{KL}_0$, which is integrable with
		$\int_0^\infty\tilde\beta(r,t)\,\mathrm dt=2r/\nu(\sqrt r)<\infty$, and
		Proposition~\ref{prop:reduction} yields the assertion.
	\end{proof}

	\section{A projection inequality}\label{app:proj-gain}

	\begin{lem}\label{lem:proj-gain}
		Let $F\subset\U$ be closed and convex, $u^\star\in F$, $\eps\in(0,1]$ and $\eta\in\U$.
		Then
		\[
		\norm{\PF(u^\star-\eps\eta)-u^\star}_\U\ge\eps\,\norm{\PF(u^\star-\eta)-u^\star}_\U .
		\]
	\end{lem}

	\begin{pf}
		Let $a\coloneqq\PF(u^\star-\eps\eta)-u^\star$ and $b\coloneqq\PF(u^\star-\eta)-u^\star$.
		The variational characterization of the projection yields
		$\ip{-\eps\eta-a}{b-a}_\U\le0$ and $\ip{-\eta-b}{a-b}_\U\le0$. Multiplying the
		second inequality by $\eps$ and adding both gives $\ip{\eps b-a}{b-a}_\U\le0$,
		hence $\eps\norm{b}^2+\norm{a}^2\le(1+\eps)\norm{a}\norm{b}$, that is,
		$(\eps\norm{b}-\norm{a})(\norm{b}-\norm{a})\le0$. If $\norm{b}>\norm{a}$, the
		first factor is nonpositive, so $\norm{a}\ge\eps\norm{b}$. Otherwise
		$\norm{a}\ge\norm{b}\ge\eps\norm{b}$.
	\end{pf}

\end{document}